\documentclass[11pt]{article}

\usepackage{amssymb}

\newcommand{\newc}{\newcommand}

\newc{\eqnoset}{\setcounter{equation}{0}}

\newcommand{\mref}[1]{(\ref{#1})}
\newcommand{\reflemm}[1]{Lemma~\ref{#1}}
\newcommand{\refrem}[1]{Remark~\ref{#1}}
\newcommand{\reftheo}[1]{Theorem~\ref{#1}}

\newcommand{\refcoro}[1]{Corollary~\ref{#1}}

\newcommand{\refsec}[1]{Section~\ref{#1}}

\newcommand{\beq}{\begin{equation}}
\newcommand{\eeq}{\end{equation}}
\newcommand{\beqno}[1]{\begin{equation}\label{#1}}

\newcommand{\barr}{\begin{array}}
\newcommand{\earr}{\end{array}}

\newc{\bearr}{\begin{eqnarray*}}
\newc{\eearr}{\end{eqnarray*}}

\newc{\bearrno}[1]{\begin{eqnarray}\label{#1}}
\newc{\eearrno}{\end{eqnarray}}

\newc{\non}{\nonumber}
\newc{\nol}{\nonumber\nl}

\newcommand{\bdes}{\begin{description}}
\newcommand{\edes}{\end{description}}
\newc{\benu}{\begin{enumerate}}
\newc{\eenu}{\end{enumerate}}
\newc{\btab}{\begin{tabular}}
\newc{\etab}{\end{tabular}}

\newtheorem{theorem}{Theorem}[section]
\newtheorem{defi}[theorem]{Definition}
\newtheorem{lemma}[theorem]{Lemma}
\newtheorem{rem}[theorem]{Remark}
\newtheorem{exam}[theorem]{Example}
\newtheorem{propo}[theorem]{Proposition}
\newtheorem{corol}[theorem]{Corollary}

\newcommand{\btheo}[1]{\begin{theorem}\label{#1}}
\newc{\brem}[1]{\begin{rem}\label{#1}\em}
\newc{\bexam}[1]{\begin{exam}\label{#1}\em}
\newc{\bdefi}[1]{\begin{defi}\label{#1}}
\newcommand{\blemm}[1]{\begin{lemma}\label{#1}}
\newcommand{\bprop}[1]{\begin{propo}\label{#1}}
\newcommand{\bcoro}[1]{\begin{corol}\label{#1}}
\newcommand{\etheo}{\end{theorem}}
\newcommand{\elemm}{\end{lemma}}
\newcommand{\eprop}{\end{propo}}
\newcommand{\ecoro}{\end{corol}}
\newc{\erem}{\end{rem}}
\newc{\eexam}{\end{exam}}
\newc{\edefi}{\end{defi}}

\newc{\rmk}[1]{{\bf REMARK #1: }}
\newc{\DN}[1]{{\bf DEFINITION #1: }}

\newcommand{\bproof}{{\bf Proof:~~}}
\newc{\eproof}{{\vrule height8pt width5pt depth0pt}\vspace{3mm}}

\newc{\bfrac}[2]{\dspl{\frac{#1}{#2}}}

\newc{\nid}{\noindent}

\newcommand{\dspl}{\displaystyle}
\newc{\grad}{\nabla}
\newc{\Div}{\mbox{div}}
\newc{\pdt}[1]{\dspl{\frac{\partial{#1}}{\partial t}}}
\newc{\pdn}[1]{\dspl{\frac{\partial{#1}}{\partial \nu}}}
\newc{\pdNi}[1]{\dspl{\frac{\partial{#1}}{\partial \mathcal{N}_i}}}
\newc{\pD}[2]{\dspl{\frac{\partial{#1}}{\partial #2}}}
\newc{\dt}{\dspl{\frac{d}{dt}}}
\newc{\bdry}[1]{\mbox{$\partial #1$}}
\newc{\sgn}{\mbox{sign}}

\newc{\Hess}[1]{\frac{\partial^2 #1}{\pdh z_i \pdh z_j}}
\newc{\hess}[1]{\partial^2 #1/\pdh z_i \pdh z_j}

\newc{\ag}{\alpha}
\newc{\bg}{\beta}
\newc{\cg}{\gamma}\newc{\Cg}{\Gamma}
\newc{\dg}{\delta}\newc{\Dg}{\Delta}
\newc{\eg}{\varepsilon}
\newc{\zg}{\zeta}
\newc{\thg}{\theta}
\newc{\llg}{\lambda}\newc{\LLg}{\Lambda}
\newc{\kg}{\kappa}
\newc{\rg}{\rho}
\newc{\sg}{\sigma}\newc{\Sg}{\Sigma}
\newc{\tg}{\tau}
\newc{\fg}{\phi}\newc{\Fg}{\Phi}
\newc{\vfg}{\varphi}
\newc{\og}{\omega}\newc{\Og}{\Omega}
\newc{\pdh}{\partial}

\newc{\ccG}{{\cal G}}

\newc{\ii}[1]{\int_{#1}}
\newc{\iidx}[2]{{\dspl\int_{#1}~#2~dx}}
\newc{\bii}[1]{{\dspl \ii{#1} }}
\newc{\biii}[2]{{\dspl \iii{#1}{#2} }}
\newc{\su}[2]{\sum_{#1}^{#2}}
\newc{\bsu}[2]{{\dspl \su{#1}{#2} }}

\newc{\biiom}[1]{{\dspl\int_{\bdrom}~ #1 ~d\sg}}
\newc{\io}[1]{{\dspl\int_{\Og}~ #1 ~dx}}
\newc{\bio}[1]{{\dspl\int_{\bdrom}~ #1 ~d\sg}}
\newc{\bsir}{\bsu{i=1}{r}}
\newc{\bsim}{\bsu{i=1}{m}}

\newc{\iibr}[2]{\iidx{\bprw{#1}}{#2}}
\newc{\Intbr}[1]{\iibr{R}{#1}}
\newc{\intbr}[1]{\iibr{\rg}{#1}}
\newc{\intt}[3]{\int_{#1}^{#2}\int_\Og~#3~dxdt}
\newc{\itQ}[2]{\dspl{\int\hspace{-2.5mm}\int_{#1}~#2~dz}}
\newc{\mitQ}[2]{\dspl{\rule[1mm]{4mm}{.3mm}\hspace{-5.3mm}\int\hspace{-2.5mm}\int_{#1}~#2~dz}}
\newc{\mitQQ}[3]{\dspl{\rule[1mm]{4mm}{.3mm}\hspace{-5.3mm}\int\hspace{-2.5mm}\int_{#1}~#2~#3}}

\newc{\mitx}[2]{\dspl{\rule[1mm]{3mm}{.3mm}\hspace{-4mm}\int_{#1}~#2~dx}}
\newc{\mitmu}[2]{\dspl{\rule[1mm]{3mm}{.3mm}\hspace{-4mm}\int_{#1}~#2~d\mu}}
\newc{\iidmu}[2]{{\dspl\int_{#1}~#2~d\mu}}

\newc{\iidm}[3]{{\dspl\int_{#1}~#2~d #3}}

\newc{\itQmu}[2]{\dspl{\int\hspace{-2.5mm}\int_{#1}~#2~d\mu}}
\newc{\mitQmu}[2]{\dspl{\rule[1mm]{4mm}{.3mm}\hspace{-5.3mm}\int\hspace{-2.5mm}\int_{#1}~#2~d\mu}}

\newc{\mitQq}[2]{\dspl{\rule[1mm]{4mm}{.3mm}\hspace{-5.3mm}\int\hspace{-2.5mm}\int_{#1}~#2~d\bar{z}}}
\newc{\itQq}[2]{\dspl{\int\hspace{-2.5mm}\int_{#1}~#2~d\bar{z}}}

\newc{\pder}[2]{\dspl{\frac{\partial #1}{\partial #2}}}

\newc{\bdrom}{\bdry{\Og}}

\newc{\bilhom}{\mbox{Bil}(\mbox{Hom}(\RR^{nm},\RR^{nm}))}
\newc{\VV}[1]{{V(Q_{#1})}}

\newc{\ccA}{{\mathcal A}}
\newc{\ccB}{{\mathcal B}}
\newc{\ccC}{{\mathcal C}}
\newc{\ccD}{{\mathcal D}}
\newc{\ccE}{{\mathcal E}}
\newc{\ccH}{\mathcal{H}}
\newc{\ccF}{\mathcal{F}}
\newc{\ccI}{{\mathcal I}}
\newc{\ccJ}{{\mathcal J}}
\newc{\ccK}{{\mathcal K}}
\newc{\ccP}{{\mathcal P}}
\newc{\ccQ}{{\mathcal Q}}
\newc{\ccR}{{\mathcal R}}
\newc{\ccS}{{\mathcal S}}
\newc{\ccT}{{\mathcal T}}
\newc{\ccX}{{\mathcal X}}
\newc{\ccY}{{\mathcal Y}}
\newc{\ccZ}{{\mathcal Z}}

\newc{\bb}[1]{{\mathbf #1}}

\newc{\myprod}[1]{\langle #1 \rangle}
\newc{\mypar}[1]{\left( #1 \right)}
\newc{\BLLg}{\mathbf{\LLg}}

\newc{\mA}{\mathbf{A}}
\newc{\mB}{\mathbf{B}}
\newc{\mC}{\mathbf{C}}
\newc{\mD}{\mathbf{D}}
\newc{\mE}{\mathbf{E}}
\newc{\mF}{\mathbf{F}}
\newc{\mJ}{\mathbf{J}}
\newc{\mG}{\mathbf{G}}
\newc{\mP}{\mathbf{P}}
\newc{\mR}{\mathbf{R}}
\newc{\mQ}{\mathbf{Q}}
\newc{\mX}{\mathbf{X}}
\newc{\muu}{\mathbf{u}}
\newc{\mvv}{\mathbf{v}}

\newc{\mllg}{\mathbb{\lambda}}
\newc{\mLLg}{\mathbf{\LLg}}

\newc{\lspn}[2]{\mbox{$\| #1\|_{\Lsp{#2}}$}}
\newc{\Lpn}[2]{\mbox{$\| #1\|_{#2}$}}
\newc{\Hn}[1]{\mbox{$\| #1\|_{H^1(\Og)}$}}

\newc{\mynorm}[2]{\| #1\|_{#2}}

\newcommand{\RR}{{\rm I\kern -1.6pt{\rm R}}}

\newc{\itQQ}[2]{\dspl{\int_{#1}#2\,dz}}
\newc{\mmitQQ}[2]{\dspl{\rule[1mm]{4mm}{.3mm}\hspace{-4.3mm}\int_{#1}~#2~dz}}
\newc{\MmitQQ}[2]{\dspl{\rule[1mm]{4mm}{.3mm}\hspace{-4.3mm}\int_{#1}~#2~d\mu}}

\newc{\MUmitQQ}[3]{\dspl{\rule[1mm]{4mm}{.3mm}\hspace{-4.3mm}\int_{#1}~#2~d#3}}
\newc{\MUitQQ}[3]{\dspl{\int_{#1}~#2~d#3}}

\newc{\mccP}{\mathbb{P}}
\newc{\mccK}{\mathbb{K}}

\newc{\DKTmU}{\mccK(U)}
\newc{\DKTmUold}{(K_U(U)^{-1})^T}

\newc{\myPi}{\mathbf{W}}
\newc{\myIbar}{\bar{\ccI}_1}
\newc{\myIhat}{\hat{\ccI}_1}
\newc{\myIbreve}{\breve{\ccI}_0}

\newc{\mmk}{\mathbf{k}}

\newc{\extraI}{\mathbb{I}}
\newc{\mccIe}{\ccI}
\newc{\mccIee}{\ccF}

\begin{document}

\vspace*{-.8in}
\begin{center} {\LARGE\em Existence of Strong Solutions to Degenerate or Singular Strongly Coupled Elliptic Systems.}

 \end{center}

\vspace{.1in}

\begin{center}

{\sc Dung Le}{\footnote {Department of Mathematics, University of
Texas at San
Antonio, One UTSA Circle, San Antonio, TX 78249. {\tt Email: Dung.Le@utsa.edu}
}}

\end{center}

\begin{abstract}
A general class of strongly coupled elliptic systems with quadratic growth in gradients is considered and the existence of their strong solutions is established. The results greatly improve those in a recent paper \cite{dleJFA} as the systems can be either degenerate or singular when their solutions become unbounded. A unified proof for both cases is presented. Most importantly, the VMO assumption in \cite{dleJFA} will be replaced by a much versatile one thanks to a new local weighted Gagliardo-Nirenberg involving BMO norms. Examples in physical models will be provided. \end{abstract}

\vspace{.2in}

\section{Introduction}\label{introsec}\eqnoset

In this paper, for any bounded domain $\Og$ with smooth boundary in $\RR^n$, $n\ge2$, we consider the following elliptic system of $m$ equations ($m\ge2$) 
\beqno{e1}\left\{\barr{l} -\Div(A(x,u)Du)=\hat{f}(x,u,Du),\quad x\in \Og,\\\mbox{$u=0$ or $\frac{\partial u}{\partial \nu}=0$ on $\partial \Og$}. \earr\right.\eeq
where $A(x,u)$ is a $m\times m$ matrix in $x\in\Og$ and $u\in\RR^m$, $\hat{f}:\Og\times\RR^m\times \RR^{mn}\to\RR^m$ is a vector valued function. 

We say that $u$ is a strong solution if $u$ is continuous on $\bar{\Og}$ with $Du\in L^\infty_{loc}(\Og)$ and $D^2u\in L^2_{loc}(\Og)$. Hence, $u$ solves \mref{e1} almost everywhere in $\Og$.

The strongly coupled system \mref{e1} appears in many applications, ranging from differential geometry to physical models. For instance, it describes the steady states of Maxwell-Stephan systems describing the diffusive transport of multicomponent mixtures, models in reaction and diffusion in electrolysis and diffusion of polymers, or population dynamics, among others. 

It is always assumed that the matrix $A(x,u)$ is elliptic in the sense that there exist two scalar positive continuous functions $\llg_1(x,u), \llg_2(x,u)$ such that \beqno{genelliptic} \llg_1(x,u)|\zeta|^2 \le \myprod{A(x,u)\zeta,\zeta}\le \llg_2(x,u)|\zeta|^2 \quad \mbox{for all } x\in\Og,\,u\in\RR^m,\,\zeta \in\RR^{nm}.\eeq If there exist positive constants $c_1,c_2$ such that $c_1\le \llg_1(x,u)$ and  $\llg_2(x,u)\le c_2$ then we say that $A(x,u)$ is {\em regular elliptic}. If $c_1\le \llg_1(x,u)$ and  $\llg_2(x,u)/\llg_1(x,u)\le c_2$, we say that $A(x,u)$ is {\em uniform elliptic}. On the other hand, if we allow $c_1=0$ and $\llg_1(x,u)$ tend to zero (respectively, $\infty$) when $|u|\to\infty$ then we say that $A(x,u)$ is {\em singular} (respectively, {\em degenerate}). 

The first fundamental problem in the study of \mref{e1} is the existence and regularity of its solutions. One can decide to work with either weak or strong solutions. In the first case, the existence of a weak solution can be achieved via Galerkin or variational methods \cite{Gius} but its regularity (e.g., boundedness, H\"older continuity of the solution and its higher derivatives) is an open issue and difficult to address. Several works (see \cite{Gius} and the reference therein) have been done along this line to establish only {\em partial regularity} of {\em bounded} weak solutions, wherever they are VMO. The assumption on the boundedness of weak solutions is a very severe and hard to check one, as maximum principles do not generally exist for systems (i.e. $m>1$) like \mref{e1}. One usually needs to use ad hoc techniques on the case by case basis to show that $u$ is bounded.  Even for bounded weak solutions, we only know that they are partially regular, i.e. H\"older continuous almost everywhere. Techniques in this approach rely heavily on the fact that $A(x,u)$ is {\em regular elliptic}, and hence the boundedness of weak solutions.

In our recent work \cite{dleJFA}, we choose a different approach making use of fixed point theory and discussing the existence of {\em strong} solutions of \mref{e1} under the weakest assumption that they are a-priori VMO, not necessarily bounded, and general structural conditions on the data of \mref{e1} which are independent of $x$, we assumed only that $A(u)$ is {\em uniformly elliptic}.  Applications were presented in \cite{dleJFA} when $\llg_1(u)$ has a positive polynomial growth in $|u|$ and, without the boundedness assumption on the solutions, so \mref{e1} can be degenerate as $|u|\to\infty$.

In this paper, we will establish much stronger results than those in \cite{dleJFA} under much more general assumptions on the structure of \mref{e1} below. Beside the minor fact that the data can depend on $x$, we allow further that:
\begin{itemize}
	\item $A(x,u)$ can be either degenerate or singular as $|u|$ tend to infinity;
	\item  $\hat{f}(x,u,Du)$ can have a quadratic growth in $Du$.
	
\end{itemize}

Most remarkably, the key assumption in \cite{dleJFA} that $u$ is VMO will be replaced by a more versatile one in this paper: $K(u)$ is VMO for some suitable map $K:\RR^m\to\RR^m$. This allows us to consider the singular case where one may not be able to estimate the BMO norm of $u$  but that of $K(u)$ in small balls. Examples of this case in applications will be provided in \refsec{res}.

One of the key ingredients in the proof in \cite{dleJFA} is the {\em local} weighted Gagliardo-Nirenberg inequality involving BMO norm \cite[Lemma 2.4]{dleANS}. This allows us to consider VMO solutions in \cite{dleJFA}. In this paper, we make use of a new version of this inequality reported in our work \cite{dleGNnew} replacing the BMO norm of $u$ by that of $K(u)$ for some suitable map $K:\RR^m\to\RR^m$.

We consider the following structural conditions on the data of \mref{e1}.

\bdes

\item[A)] $A(x,u)$ is $C^1$ in $x\in\Og$, $u\in\RR^m$ and there exist a constant $C_*>0$  and  scalar $C^1$ positive functions $\llg(u),\og(x)$  such that  for all $u\in\RR^m$,  $\zeta\in\RR^{mn}$ and $x\in\Og$ 
\beqno{A1} \llg(u)\og(x)|\zeta|^2 \le \myprod{A(x,u)\zeta,\zeta} \mbox{ and } |A(x,u)|\le C_*\llg(u)\og(x).\eeq 

In addition, there is a constant $C$ such that $|\llg_u(u)||u|\le C\llg(u)$ and  \beqno{Axcond} |A_u(x,u)|\le C|\llg_u(u)|\og(x),\; |A_x(x,u)|\le C|\llg(u)||D\og|.\eeq

\edes

Here and throughout this paper, if $B$ is a $C^1$ (vector valued) function in $u\in \RR^m$ then we abbreviate its  derivative $\frac{\partial B}{\partial u}$ by $B_u$. Also, with a slight abuse of notations, $A(x,u)\zeta$, $\myprod{A(x,u)\zeta,\zeta}$ in \mref{genelliptic}, \mref{A1} should be understood in the following way: For $A(x,u)=[a_{ij}(x,u)]$, $\zeta\in\RR^{mn}$ we write $\zeta=[\zeta_i]_{i=1}^m$ with $\zeta_i=(\zeta_{i,1},\ldots\zeta_{i,n})$ and $$A(x,u)\zeta=[\Sigma_{j=1}^m a_{ij}\zeta_j]_{i=1}^m,\; \myprod{A(x,u)\zeta,\zeta}=\Sigma_{i,j=1}^m a_{ij}\myprod{\zeta_i,\zeta_j}.$$

We also assume that $A(x,u)$ is regular elliptic for {\em bounded} $u$.

\bdes\item[AR)] $\og\in C^1(\Og)$ and there are positive numbers $\mu_*,\mu_{**}$ such that 
\beqno{Aregcond2} \mu_*\le\og(x)\le \mu_{**}, \; |D\og(x)|\le \mu_{**} \quad \forall x\in\Og.\eeq For any bounded set $K\subset\RR^m$ there is a constant $\llg_*(K)>0$ such that
\beqno{Aregcond1} \llg_*(K)\le\llg(u) \quad \forall u\in K.\eeq \edes

Concerning the reaction term $\hat{f}(x,u,Du)$, which may have linear or {\em quadratic} growth in $Du$, we assume the following condition.
\bdes \item[F)] There exist a constant $C$ and a nonegative differentiable function $f:\RR^m\to\RR$ such that $\hat{f}$ satisfies: For  any diffrentiable vector valued functions $u:\RR^n\to\RR^m$ and $p:\RR^n\to\RR^{mn}$ we assume either that
\bdes\item[f.1)] $\hat{f}$ has a linear growth in $p$ 
\beqno{FUDU11}|\hat{f}(x,u,p)| \le C\llg(u)|p|\og(x) + f(u)\og(x),\eeq $$|D\hat{f}(x,u,p)| \le  C(\llg(u)|Dp|+ |\llg_u(u)||p|^2)\og+ C\llg(u)|p||D\og|+C|D(f(u)\og(x))|;$$\edes or \bdes 
\item[f.2)] $\llg_{uu}(u)$ exists and $\hat{f}$ has a quadratic growth in $p$ \beqno{FUDU112}|\hat{f}(x,u,p)| \le C|\llg_u(u)||p|^2\og(x) + f(u)\og(x),\eeq  $$\barr{lll}|D\hat{f}(x,u,p)| &\le&  C(|\llg_u(u)||p||Dp|+ |\llg_{uu}(u)||p|^3)\og+ C|\llg_u(u)||p|^2|D\og|\\&&+C|D(f(u)\og(x))|.\earr$$ Furthermore, we assume that \beqno{llgquadcond} |\llg_{uu}(u)|\llg(u)\le C|\llg_u(u)|^2.\eeq
\edes
\edes

By a formal differentiation of \mref{FUDU11} and \mref{FUDU112}, one can see that  the growth conditions for $\hat{f}$ naturally implies those of  $D\hat{f}$ in the above assumption. The condition \mref{llgquadcond} is verified easily if $\llg(u)$ has a polynomial growth in $|u|$.

We organize our paper as follows. In \refsec{res} we state our main results and their applications which are actually consequences of the most general but technical \reftheo{gentheo1} in \refsec{w12est} where we deal with general map $K$. The proof of the results in \refsec{res}, which provide examples of the map $K$, will thus be provided in \refsec{proofmainsec}. In \refsec{GNsec} we state the new version of the local weighted Gagliardo-Nirenberg inequality in \cite{dleGNnew} to prepare for the proof the main technical theorem in \refsec{w12est}. We collect some elementary but useful facts in our proof in \refsec{appsec}.

\section{Preliminaries and Main Results}\eqnoset\label{res}

We state the main results of this paper in this section. In fact, these results are consequences of our main technical results in \refsec{w12est} assuming general conditions A) and F) and, roughly speaking, some a~priori knowledge on the smallness of the BMO norm of $K(u)$ for a general map $K:\RR^m\to\RR^m$ and any strong solution $u$ to \mref{e1}.

We define the measure $d\mu =\og(x)dx$ and recall that a vector valued function $f\in L^1(\Og,\mu)$ is said to be in $BMO(\Og,\mu)$ if \beqno{bmodef} [f]_{*,\mu}:=\sup_{B_R\subset\Og}\mitmu{B_R}{|f-f_{B_R}|}<\infty,\quad f_{B_R}:=\frac{1}{\mu(B_R)}\iidmu{B_R}{f}.\eeq We then define $$\|f\|_{BMO(\Og,\mu)}:=[f]_{*,\mu}+\|f\|_{L^1(\Og,\mu)}.$$

Throughout this paper, in  our statements and proofs, we use $C,C_1,\ldots$ to denote various constants which can change from line to line but depend only on the parameters of the hypotheses in an obvious way. We will write $C(a,b,\ldots)$ when the dependence of  a constant $C$ on its parameters is needed to emphasize that $C$ is bounded in terms of its parameters. We also write $a\lesssim b$ if there is a universal constant $C$ such that $a\le Cb$. In the same way, $a\sim b$ means $a\lesssim b$ and $b\lesssim a$.

To begin, as in \cite{dleANS} with $A$ is independent of $x$, we assume that the eigenvalues of the matrix $A(x,u)$ are not too far apart. Namely, for $C_*$ defined in \mref{A1} of A) we assume
\bdes\item[SG)] $(n-2)/n <C_*^{-1}$. \edes Here $C_*$ is, in certain sense, the ratio of the largest and smallest eigenvalues of $A(x,u)$. This condition seems to be necessary as we deal with systems, cf.  \cite{dleN}.

In this section we assume further the following growth conditions on $\llg(u)$ and $f(u)$. \bdes\item[AF)] There is a constant $k$ such that for all $u\in\RR^m$
\beqno{llucond} \llg(u)\sim (\llg_0+|u|)^k \mbox{ and }|\llg_u(u)|\sim (\llg_0+|u|)^{k-1}. \eeq

If $\hat{f}$ depends on $x$, then we also assume that there is a constant $C$ such that
\beqno{ffucond} |f(u)|\le C|f_u(u)|(1+|u|).\eeq
\edes

We imbed \mref{e1} in the following family of systems 
\beqno{e1famzzz}\left\{\barr{l} -\Div(A(x,\sg u)Du)=\hat{f}(x,\sg u,\sg Du),\quad x\in \Og, \sg\in[0,1],\\\mbox{$u=0$ or $\frac{\partial u}{\partial \nu}=0$ on $\partial \Og$}. \earr\right.\eeq

Our first main result is the following.

\btheo{dleNL-mainthm} 
Assume that A), AR), F), AF) and SG) hold. Suppose also that the following integrability conditions hold for {\em any strong} solution $u$ of \mref{e1famzzz}: There exist $r_0>n/2$, $\bg_0\in(0,1)$    such that the following quantities
\beqno{llgmainhyp009} \|\llg^{-1}(u)\|_{L^{\frac n2}(\Og,\mu)},\;\||f_u(u)|\llg^{-1}(u)\|_{L^{r_0}(\Og,\mu)},\; \|(\llg(u)|u|^2)^{\bg_0}\|_{L^{1}(\Og,\mu)},\eeq 
\beqno{llgmainhyp0099}\iidmu{\Og}{(|f_u(u)|+\llg(u))|Du|^{2}}\eeq are bounded by some constant $C_0$.

Define $K_0:\RR^m\to\RR^m$ by $K_0(u)=|\log(|U|)||U|^{-1}U$, $U=[\llg_0+|u_i|]_{i=1}^m$. We assume that the BMO norms of $K_0(u)$ and  $\log(\llg_0+|u|)$ are small in small balls in the following sense: For any $\eg>0$ there is $R_\eg>0$ such that  \beqno{mainloghyp} \|\log(\llg_0+|u|)\|_{BMO(B_R,\mu)},\;\|K_0(u)\|_{BMO(B_R,\mu)}<\eg \quad \mbox{for all $B_R\subset\Og$ with $R\le R_\eg$}.\eeq 

Then \mref{e1} has a strong solution $u$.

 \etheo

The condition \mref{mainloghyp} on the smallness of the BMO norm of $K_0(u)$ in small balls is the most crucial one in applications. The next result is more applicable in the checking of this condition. 

\bcoro{dlec1coro} The conclusion of \reftheo{dleNL-mainthm} holds if \mref{mainloghyp} is replaced by: There exist $\ag\in(0,1)$ and a constant $C_0$ such that 
\beqno{pis20az5} \iidmu{\Og}{(\llg_0+|u|)^{-n\ag}|Du|^{n}}\le C_0.\eeq

\ecoro

In \cite{dleJFA}, we consider the case $\llg(u)\sim(\llg_0+|u|)^k$ with $k>0$ and assume that $u$ has small BMO norm in small balls, which can be verified by establishing that $\|Du\|_{L^n(\Og)}$ is bounded. The assumption \mref{pis20az5} is of course much weaker, especially when $|u|$ is large, and can apply to the case $k<0$.

We present an application of \refcoro{dlec1coro}. This example concerns cross diffusion systems with polynomial growth data on planar domains. This type of systems occurs in many applications in mathematical biology and ecology. We will see that the assumptions of the corollary can be verified by a very simple integrability assumption on the solutions.

\bcoro{2dthm1} Let $n=2$. Suppose A), F) and $f(u)\lesssim (\llg_0+|u|)^l$ and $\llg(u)\sim(\llg_0+|u|)^{k}$ for some $k,l$ satisfying \beqno{klcond} k>\frac{-2C_*}{C_*-1} \mbox{ and } l-k<\frac{C_*+1}{C_*-1}.\eeq 
If $\hat{f}$ has a quadratic growth in $Du$ as in \mref{FUDU112} of f.2) then we assume further that \beqno{FUDU112z} |\hat{f}(x,u,p)| \le \eg_0|\llg_u(u)||p|^2 + f(u),\eeq with $\eg_0$ being sufficiently small.

If there is a constant $C_0$ such that for {\em any strong} solution $u$ of \mref{e1famzzz} \beqno{keyn2lnorm0} \|u\|_{L^{l_0}(\Og,\mu)}\le C_0\quad \mbox{for some $l_0>\max\{l,l-k-1\}$},\eeq then \mref{e1} has a strong solution $u$.  \ecoro

The assumption \mref{keyn2lnorm0} is a very weak one. For example, if $k\ge-1$ then we see from the growth condition $|f(u)|\lesssim (\llg_0+|u|)^l$ that  \mref{keyn2lnorm0} simply requires that $l_0>l$, or equivalently, $f(u)\in L^r(\Og,\mu)$ for some $r>1$. 

This result greatly generalizes \cite[Corollary 3.9]{dleJFA} in many aspects. Beside the fact that we allow quadratic growth in $Du$ for $\hat{f}(x,u,Du)$ and $k<0$, we also consider a much general relation between the the growths of $f(u)$ and $\llg(u)$ in \mref{klcond}, while we assume in \cite{dleJFA} that $f(u)\lesssim \llg(u)|u|$ (i.e. $l-k=1$).

In the second main result, we consider the following generalized SKT system (see \cite{dleJFA,SKT,yag}) with Dirichlet or Neumann boundary conditions on a bounded domain $\Og\subset\RR^n$ with $n\le4$.
\beqno{genSKT} -\Delta(P_i(u))=B_i(u,Du)+ f_i(u), \quad i=1,\ldots,m.\eeq Here, $P_i:\RR^m\to\RR$ are $C^2$ functions. The functions $B_i, f_i$ are $C^1$ functions on $\RR^m\times\RR^{mn}$ and $\RR^m$ respectively. We will assume that $B_i(u,Du)$ has linear growth in $Du$.

By a different choice of the map $K$ in the main technical theorem, we have the following.

\btheo{n3SKT} Assume that the matrix $A(u)=\mbox{diag}[(P_i)_u(u)]$ satisfies the condition A) with $\llg(u)\sim (\llg_0+|u|)^k$ for some $k\ge-1$. Moreover, $\hat{f}(u,Du)=\mbox{diag}[B_i(u,Du)+f_i(u)]$ satisfies the following special version f.1) of F)  $$|B_i(u,Du)|\le C\llg(u)|Du|, \; |f_i(u)|\le f(u).$$

Thus, \mref{genSKT} can be written as \mref{e1}. Assume that there exist $r_0>n/2$ and a constant $C_0$  such that for {\em any strong} solution $u$ of \mref{e1famzzz}
\beqno{fuhyp0}\|f_u(u)\llg^{-1}(u)\|_{L^{r_0}(\Og)}\le C_0,\eeq and the following conditions. \bdes\item[i)] If $k\ge0$ then  $\|u\|_{L^1(\Og)}\le C_0$. \item[ii)] If $k\in[-1,0)$ then $\|u\|_{L^{-kn/2}(\Og)}\le C_0$. Furthermore,
\beqno{fuhyp0k}\|\llg^{-2}(u)f_u(u)\|_{L^\frac{n}{2}(\Og)}\le C_0.\eeq
\edes

Then 
\mref{genSKT} has a strong solution for $n=2,3,4$.
\etheo

The above result generalizes \cite[Corollary 3.10]{dleJFA} where we assumed that $\hat{f}$ is independent of $Du$, $k>0$ and $f(u)\lesssim \llg(u)|u|$. In this case, it is natural to assume that $|f_u(u)|\lesssim\llg(u)$ so that \mref{fuhyp0} obviously holds. We also have $|\llg^{-2}(u)f_u(u)|\lesssim \llg^{-1}(u)$ so that \mref{fuhyp0k} is in fact a consequence of the assumption $\|u\|_{L^{-kn/2}(\Og)}\le C_0$ in ii).

We should remark that all the assumptions on strong solutions of the family \mref{e1famzzz} can be checked by considering the case $\sg=1$ (i.e. \mref{e1} or \mref{genSKT}) because these systems satisfy the same structural conditions uniformly with respect to the parameter $\sg\in[0,1]$.

\section{A general  local weighted Gagliardo-Nirenberg inequality} \eqnoset\label{GNsec}
In this section, we present a  local weighted Gagliardo-Nirenberg inequality in our recent work \cite{dleGNnew}, which will be one of the main ingredients of the proof of our main technical theorem in \refsec{w12est}. 
This inequality generalizes \cite[Lemma 2.4]{dleANS} by replacing the Lebesgue measure with general one and the BMO norm of $u$ with that of $K(u)$ where $K$ is a suitable map on $\RR^m$, and so the applications of our main technical theorem in the next section will be much more versatile than those in \cite{dleJFA,dleANS}.

Let us begin by describing the assumptions in \cite{dleGNnew} for this general inequality.
We need to recall some well known notions from Harmonic Analysis.

Let $\og\in L^1(\Og)$ be a nonnegative function and  define the measure $d\mu=\og(x)dx$. For any $\mu$-measurable subset $A$ of $\Og$  and any  locally $\mu$-integrable function $U:\Og\to\RR^m$ we denote by  $\mu(A)$ the measure of $A$ and $U_A$ the average of $U$ over $A$. That is, $$U_A=\mitmu{A}{U(x)} =\frac{1}{\mu(A)}\iidmu{A}{U(x)}.$$

We say that $\Og$ and $\mu$ support a $q_*$-Poincar\'e inequality if the following holds. \bdes \item[P)]  There exist $q_*\in(0,2]$, $\tau_*\ge 1$ and some constant $C_P$ such that \beqno{Pineq2} 
\mitmu{B}{|h-h_{B}|}\le C_Pl(B)
\left(\mitmu{\tau_*B}{|Dh|^{q_*}}\right)^\frac{1}{q_*}\eeq
for  any cube $B\subset\Og$ with side length $l(B)$ and any function $u\in C^1(B)$. \edes

Here and throughout this section,  we denote by $l(B)$ the side length of $B$ and by $\tau B$ the cube which is concentric with $B$ and has side length $\tau l(B)$. We also write $B_R(x)$ for a cube centered at $x$ with side length $R$ and sides parallel to to standard axes of $\RR^n$. We will omit $x$ in the notation $B_R(x)$ if no ambiguity can arise.

We consider the following conditions on the density $\og(x)$.

\bdes \item[LM.1)]  For some $N\in(0,n]$ and any ball $B_r$ we have $\mu(B_r)\le C_\mu r^N$. 
Assume also that $\mu$ supports the 2-Poincar\'e inequality \mref{Pineq2} in P). Furthermore, $\mu$ is doubling and satisfies  the following inequality for some $s_*>0$ \beqno{fracmu} \left(\frac{r}{r_0}\right)^{s_*}\le C_\mu\frac{\mu(B_r(x))}{\mu(B_{r_0}(x_0))},\eeq where $B_r(x), B_{r_0}(x_0)$ are any cubes with $x\in B_{r_0}(x_0)$.

\item[LM.2)]  $\og=\og_0^2$ for some $\og_0\in C^1(\Og)$ and $d\mu=\og_0^2 dx$  also supports a Hardy type inequality: There is a constant $C_H$ such that for any function $u\in C^1_0(B)$\beqno{lehr1m} \iidx{\Og}{|u|^2|D\og_0|^2}\le C_H\iidx{\Og}{|Du|^2\og_0^2}.\eeq 
\edes

For $\cg\in(1,\infty)$ we say that a nonnegative locally integrable function $w$ belongs to the class $A_\cg$ or $w$ is an $A_\cg$ weight on $\Og$ if the quantity
\beqno{aweight} [w]_{\cg,\Og} := \sup_{B\subset\Og} \left(\mitmu{B}{w}\right) \left(\mitmu{B}{w^{1-\cg'}}\right)^{\cg-1} \quad\mbox{is finite}.\eeq
Here, $\cg'=\cg/(\cg-1)$. For more details on these classes we refer the reader to \cite{OP,st}. If the domain $\Og$ is specified we simply denote $ [w]_{\cg,\Og}$ by $ [w]_{\cg}$.

We assume the following hypotheses.

\bdes
\item[A.1)] Let $K:\mbox{dom}(K)\to\RR^m$ be a $C^1$ map on a domain $\mbox{dom}(K)\subset\RR^m$ such that $\DKTmU=(K_U(U)^{-1})^T$ exists and $\mccK_U\in L^\infty(\mbox{dom}(K))$.

Furthermore, let $\Fg,\LLg:\mbox{dom}(K)\to\RR^+$ be $C^1$ positive functions. We assume that for all $U\in \mbox{dom}(K)$
\beqno{kappamain} |\DKTmU|\lesssim \LLg(U)\Fg^{-1}(U),\eeq
\beqno{logcondszmain}|\Fg_U(U)||\mccK(U)|\lesssim \Fg(U).\eeq
\edes

Let $\Og_*$ be a proper subset of $\Og$ and $\og_*$ be a function in $C^1(\Og)$ satisfying  
\beqno{subogm} \og_*\equiv 1 \mbox{ in $\Og_*$ and } \og_*\le 1 \mbox{ in $\Og$}.\eeq

For any $U\in C^2(\Og,\mbox{dom}(K))$  we denote
\beqno{Idefm} I_1:=\iidmu{\Og}{\Fg^2(U)|DU|^{2p+2}},\;
I_2:=\iidmu{\Og}{\LLg^2(U)|DU|^{2p-2}|D^2U|^2},\eeq
\beqno{Idef1zm} \myIbar:=\iidmu{\Og}{|\LLg_U(U)|^2|DU|^{2p+2}},\;I_{1,*}:=\iidmu{\Og_*}{\Fg^2(U)|DU|^{2p+2}},\eeq \beqno{I0*m} \breve{I}_{0,*}:=\sup_\Og|D\og_*|^2\iidmu{\Og}{\LLg^2(U)|DU|^{2p}}.\eeq 

By \refrem{PSrem} below, the assumption PS) in \cite{dleGNnew} that $\mu$ supports a Poincar\'e-Sobolev inequality  is then satisfied. We established the following local weighted Gagliardo-Nirenberg inequality in \cite{dleGNnew}.

\btheo{GNlocalog1m} Suppose LM.1)-LM.2), A.1). Let $U\in C^2(\Og,\mbox{dom}(K))$ and satisfy
\beqno{boundaryzm}\myprod{\og_*\og_0^2 \Fg^2(U)\DKTmU DU,\vec{\nu}}=0\eeq on $\partial\Og$ where $\vec{\nu}$ is the outward normal vector of $\partial\Og$.  Let $\myPi(x):=\LLg^{p+1}(U(x))\Fg^{-p}(U(x))$ and assume that $[\myPi^{\ag}]_{\bg+1}$ is finite for some $\ag>2/(p+2)$ and $\bg<p/(p+2)$.

Then, for any $\eg>0$ there are constants $C, C([\myPi^{\ag}]_{\bg+1})$ such that
\beqno{GNlocog11m}I_{1,*}\le  \eg I_1+\eg^{-1}C\|K(U)\|_{BMO(\mu)}^2[I_2+\myIbar+C([\myPi^{\ag}]_{\bg+1}) [I_2+\myIbar+\breve{I}_{0,*}]].\eeq
Here, $C$ also depends on $C_{P},C_\mu$ and $C_H$.

\etheo

For our purpose in this paper we need only a special case of \reftheo{GNlocalog1m} where  $\Og,\Og_*$ are concentric balls $B_s,B_t$, $0<s<t$.
We let $\og_*$ be a cutoff function for $B_s,B_{t}$:  $\og_*$ is a $C^1$ function satisfying $\og_*\equiv1$ in $B_s$ and $\og_*\equiv0$ outside $B_t$ and $|D\og_*|\le 1/(t-s)$. The condition \mref{boundaryzm} of the above theorem is clearly satisfied on the boundary of $\Og=B_t$. We also consider only the case $\Fg(U)\sim|\LLg_U(U)|$.

We then have the following corollary.

\bcoro{GNlocalog1mcoro} Suppose that LM.1)-LM.2), A.1) hold for   $\Fg(U)=|\LLg_U(U)|$. Accordingly, define  $\myPi_p(x):=\LLg^{p+1}(U(x))|\LLg_U(U(x))|^{-p}$ and let $B_t(x_0)$ be any ball in $\Og$ and assume that \bdes\item[A.2)]  $[\myPi_p^{\ag}]_{\bg+1,B_t(x_0)}$ is finite for some $\ag>2/(p+2)$ and $\bg<p/(p+2)$.\edes

 We denote (compare with \mref{Idefm}-\mref{I0*m}) \beqno{Idefmt} I_{0}(t,x_0):=\iidmu{B_t(x_0)}{\LLg^2(U)|DU|^{2p}},\; I_1(t,x_0):=\iidmu{B_t(x_0)}{|\LLg_U(U)|^2|DU|^{2p+2}},\;\eeq
\beqno{I0*mt} 
I_2(t,x_0):=\iidmu{B_t(x_0)}{\LLg^2(U)|DU|^{2p-2}|D^2U|^2}.\eeq

Then, for any $\eg>0$ and any ball $B_s(x_0)$, $0<s<t$, there are constants $C, C([\myPi_p^{\ag}]_{\bg+1,B_t(x_0)})$ with $C$ also depending on $C_{PS},C_\mu$ and $C_H$ such that for $$C_{\eg,U,\myPi}=\eg+\eg^{-1}C\|K(U)\|_{BMO(B_t(x_0),\mu)}^2[1+C([\myPi_p^{\ag}]_{\bg+1,B_t(x_0)})]$$ we have
\beqno{GNlocog11mcoro}I_{1}(s,x_0)\le  C_{\eg,U,\myPi}[I_1(t,x_0)+I_2(t,x_0)+(t-s)^{-2}I_{0}(t,x_0)].\eeq 

\ecoro

We end this section by some remarks on  the measure $\mu$.
\brem{PSrem} As we assume in LM.1) that
$\mu$ is doubling and supports a $2$-Poincar\'e inequality \mref{Pineq2}, from \cite[Section 3]{Haj} we see that  an improved $q_*$-Poincar\'e inequality also holds for any $q_*\in(p_1,2]$ for some $p_1<2$. That is, see \cite[eqn. (5)]{Haj}, there are some constants $C_P$ and $\tau_*\ge1$ such that the following $q_*$-Poincar\'e inequality holds true
\beqno{Pineq} 
\mitmu{B}{|h-h_{B}|}\le C_Pl(B)
\left(\mitmu{\tau_*B}{|Dh|^{q_*}}\right)^\frac{1}{q_*}.\eeq 
This and the assumption \mref{fracmu} and \cite[1) of Theorem 5.1]{Haj}  show that the following Poincar\'e-Sobolev inequality holds for {\em some} $\pi_*>2$ and $q_*<2$ and some constant $C_{PS}$
\beqno{PSineq} \frac{1}{l(B)}
\left(\mitmu{B}{|u-u_{B}|^{\pi_*}}\right)^\frac1{\pi_*}  \le C_{PS}
\left(\mitmu{\tau_*B}{|Du|^{2}}\right)^\frac{1}{2},\quad \pi_*>2.\eeq\erem  In fact, if $q_*<s_*$, the exponent in \mref{fracmu}, then \cite[1) of Theorem 5.1]{Haj} establishes \mref{PSineq}
for  $\pi_*=s_*q_*/(s_*-q_*)$. Thus, $\pi_*>2$ if $s_*<2q_*/(2-q_*)$. This is the case if we choose $q_*<2$ and closed to 2. Hence, the assumption PS) in \cite{dleGNnew} that $\mu$ supports a Poincar\'e-Sobolev inequality \mref{PSineq} is then satisfied for some $q_*<2$ and $\pi_*>2$  (the dimensional parameters $d,n$ in that paper are now denoted by $n,N$ respectively).

\brem{PSremH}  If $q_*=s_*$, \cite[2) of Theorem 5.1]{Haj} shows that \mref{PSineq} holds true
for any $\pi_*>1$. In addition, if $q_*>s_*$ then  the H\"older norm of $u$ is bounded in terms of $\|Du\|_{L^{q_*}(\Og,\mu)}$. 

\erem

\section{The main technical theorem}\eqnoset\label{w12est}

In this section, we establish the main result of this paper. We consider the following system
\beqno{gensys}\left\{\barr{l} -\Div(A(x,u)Du)=\hat{f}(x,u,Du),\quad x\in \Og,\\\mbox{$u=0$ or $\frac{\partial u}{\partial \nu}=0$ on $\partial \Og$}. \earr\right.\eeq

We imbed this system in the following family of systems 
\beqno{gensysfam}\left\{\barr{l} -\Div(A(x,\sg u)Du)=\hat{f}(x,\sg u,\sg Du),\quad x\in \Og, \sg\in[0,1],\\\mbox{$u=0$ or $\frac{\partial u}{\partial \nu}=0$ on $\partial \Og$}. \earr\right.\eeq

First of all, we will assume that the system \mref{gensys} satisfies the structural conditions A) and F). Additional assumptions serving the purpose of this paper then follow  for the validity of the local weighted Gagliardo-Nirenberg inequality of \refcoro{GNlocalog1mcoro} with  $\LLg(U)=\llg^\frac12(U)$.

\bdes 
\item[H)] There is a $C^1$ map $K:\RR^m\to\RR^m$ such that $\mathbb{K}(u)=(K_u(u)^{-1})^T$ exists and $\mathbb{K}_U\in L^\infty(\RR^m)$.
Furthermore,  for all $u\in \RR^m$
\beqno{kappamainz} |\mathbb{K}(u)|\lesssim \llg(u)|\llg_u(u)|^{-1}.\eeq
\edes

\brem{polyrem} We can see that the condition H) implies the condition A.1) in \reftheo{GNlocalog1m}, and then \refcoro{GNlocalog1mcoro}
with $\LLg(u)=\llg^\frac12(u)$ and $\Fg(u)=|\LLg_u(u)|$, \mref{GNlocog11mcoro} is then applicable. Indeed, the assumption \mref{kappamain} in this case is \mref{kappamainz}.  It is not difficult to see that the assumption in f.2) that  $|\llg_{uu}(u)|\llg(u)\lesssim |\llg_u(u)|^2$ and \mref{kappamainz} imply $|\Fg_u(u)||\mathbb{K}(u)|\lesssim \Fg(u)$, which gives \mref{logcondszmain} of A.1). Hence, A.1) holds by H). In particular, if $\llg$ has a polynomial growth in $u$, i.e. $\llg(u)\sim (\llg_0+|u|)^k$ for some $k\ne0$ and $\llg_0\ge0$, then H) reduced to the simple condition $|\mathbb{K}(u)|\lesssim|u|$.
\erem

For any strong solution $u$ of \mref{gensysfam} we will consider  the following assumptions. The exponents $s_*,\pi_*$ are defined in \mref{fracmu} and  in the Poincar\'e-Sobolev inequality \mref{PSineq}.

\bdes 

\item[M.0)] There exist a constant $C_0$  and some $r_0>r_*:={\pi_*}/({\pi_*}-2)$ such that
\beqno{llgmainhyp0} \|\llg^{-1}(u)\|_{L^{r_*}(\Og,\mu)},\;\|\llg^{-1}(u)f_u(u)\|_{L^{r_0}(\Og,\mu)}\le C_0,\eeq
\beqno{llgB3}\iidmu{\Og}{(|f_u(u)|+\llg(u))|Du|^{2}}\le C_0.\eeq

\item[M.1)] For any given  
$\mu_0>0$ 
there is positive $R_{\mu_0}$ sufficiently small in terms of the constants in A) and F) such that \beqno{Keymu0} \sup_{x_0\in\bar{\Og}}\|K(u)\|_{BMO(B_{R}(x_0)\cap\Og,\mu)}^2 \le \mu_0.\eeq  

Furthermore, for  $\myPi_p(x):= \llg^{p+\frac12}(u)|\llg_u(u)|^{-p}$ and any $p\in[1,s_*/2]$ there exist some $\ag>2/(p+2)$, $\bg<p/(p+2)$  such that $[\myPi_p^{\ag}]_{\bg+1,B_{R_{\mu_0}}(x_0)\cap\Og}\le C_0$. 
 \edes

The main theorem of this section is the following.

\btheo{gentheo1} Assume A), F), AR) and H).  Moreover, if $\hat{f}$ has a quadratic growth in $Du$ as in f.2) then we assume also that $n\le3$. Suppose also that any strong solution $u$ to \mref{gensysfam} satisfies M.0), M.1) for some constant $C_0$ and \beqno{hypoiterpis1z00}\iidmu{\Og}{|f_u(u)|(1+f(u)^{s_*}|f_u(u)|^{-s_*})}\le C_0,\eeq
uniformly in $\sg$.  Then the system \mref{gensys} has a strong solution. \etheo

The condition \mref{hypoiterpis1z00} is not needed if $\hat{f}$ is independent of $x$ (see \refrem{hypoiterpis1zrem}). On the other hand,  we can assume a mild condition on the growth of $f(u)$ and show that \mref{hypoiterpis1z00} can be greatly relaxed in the following result.

\bcoro{gentheo1coro} Assume as in \reftheo{gentheo1} and in addition that 
\beqno{specfcond}f(u)\lesssim|f_u(u)|(1+|u|).\eeq Then the conclusion of the theorem still holds if we replace \mref{hypoiterpis1z00} by
\beqno{hypoiterpis1z11} \|f_u(u)\|_{L^1(\Og,\mu)}\le C_0,\eeq and assume that there exists $\bg_0\in(0,1)$ such that
\beqno{hypoiterpis1z11z}\iidmu{\Og}{(\llg(u)|u|^2)^{\bg_0}}\le C_0.\eeq
\ecoro

We only discuss the existence of strong solutions in this paper so that the condition on the regularity of $\og$ in AR) seems to be necessary. Under this condition, for $\og_0=\sqrt{\og}$, the density $\og$ clearly satisfies LM.1)-LM.2) and supports the Hardy and Poincar\'e-Sobolev inequalities \mref{lehr1m} and \mref{PSineq},  with $s_*=n$ and $\pi_*=2n/(n-2)$. However, in the proof of this section, we prefer to keep the notations $s_*,\pi_*$ as in \refsec{GNsec} because many results in this section will be applicable in our future works where we will consider systems which are degenerate and singular in $u,x$ altogether, see \refrem{murem1}.

The proof of \reftheo{gentheo1} relies on the Leray Schauder fixed point index theorem. We obtain the existence of a strong solution $u$ of \mref{gensys} as a fixed point of a nonlinear map defined on an appropriate Banach space $\mX$. The proof will be based on several lemmas and we will sketch the main steps below.

We will show in \reflemm{dleANSpropt0} that there exist  $p>s_*/2$ and a constant $M_*$ depending only on the constants in A) and F) such that any strong solution $u$ of \mref{gensysfam} will satisfy \beqno{M*def} \|Du\|_{L^{2p}(\Og,\mu)}\le M_*.\eeq  This and \refrem{PSremH} imply that there are positive constants $\ag,M_0$ such that \beqno{M0def} \|u\|_{C^{\ag}(\Og)}\le M_0.\eeq

For $\sg\in[0,1]$ and any $u\in \RR^m$ and $\zeta\in\RR^{mn}$ we  define  the vector valued functions $F^{(\sg)}$ and $f^{(\sg)}$ by \beqno{Bfdef}F^{(\sg)}(x,u,\zeta):=\int_{0}^{1}\partial_\zeta F(\sg,u,t\zeta)\,dt,\quad f^{(\sg)}(x,u):=\int_{0}^{1}\partial_u F(\sg,x,tu,0)\,dt.\eeq

 For   any given $u,w\in W^{1,2}(\Og)$ we write \beqno{Bfdefalin}\mathbf{\hat{f}}(\sg,x, u,w)=F^{(\sg)}(x,u,Du)Dw+f^{(\sg)}(x,u)w+\hat{f}(x,0,0).\eeq

For any given $u\in W^{1,2}(\Og)$ satisfying \mref{M0def} we consider the following linear systems, noting that $\mathbf{\hat{f}}(\sg,x,u,w)$ is linear in $w,Dw$
\beqno{Tmapdef}\left\{\barr{l} -\Div(A(x,\sg u)Dw)+Lw=\mathbf{\hat{f}}(\sg,x,u,w)+Lu\quad x\in \Og, \\\mbox{$w=0$ or $\frac{\partial w}{\partial \nu}=0$ on $\partial \Og$}. \earr\right.\eeq
Here, $L$ is a suitable positive definite matrix depending on the constant $M_0$ such that the above system has a unique weak solution $w$ if $u$ satisfies \mref{M0def}. We then define $T_\sg(u)=w$ and apply the Leray-Schauder fixed point theorem to establish the existence of a fixed point of $T_1$.
It is clear from \mref{Bfdefalin} that $\hat{f}(x,\sg u,\sg Du)=\mathbf{\hat{f}}(\sg,x,u,u)$.  
Therefore,  from  the definition of $T_\sg$ we see that a fixed point of $T_\sg$ is a weak solution of \mref{gensysfam}. By an appropriate choice of $\mX$, we will show that these fixed points are strong solutions of \mref{gensysfam}, and so a fixed point of $T_1$ is a strong solution of \mref{gensys}.

From the proof of Leray-Schauder fixed point theorem in \cite[Theorem 11.3]{GT}, we need to  find some ball $B_M$ of radius $M$ and centered at $0$ of $\mX$ such that $T_\sg: \bar{B}_M\to \mX$ is compact and  that $T_\sg$ has no fixed point on the boundary of $B_M$.  The topological degree $\mbox{ind}(T_\sg, B_M)$ is then well defined and invariant by homotopy so that $\mbox{ind}(T_1, B_M)=\mbox{ind}(T_0, B_M)$. It is easy to see that the latter is nonzero because the linear system $$\left\{\barr{l} -\Div(A(x,0)Du)=\mathbf{\hat{f}}(x,0,0)\quad x\in \Og, \\\mbox{$u=0$ or $\frac{\partial u}{\partial \nu}=0$ on $\partial \Og$}, \earr\right.$$ has a unique solution in $B_M$. Hence, $T_1$ has a fixed point in $B_M$.

Therefore, the theorem is proved as we will establish the following claims.

\bdes \item[Claim 1] There exist a Banach space $\mX$ and $M>0$ such that the map $T_\sg:\bar{B}_M\to\mX$ is well defined and compact. 
\item[Claim 2] $T_\sg$ has no fixed point on the boundary of $\bar{B}_M$. That is,  $\|u\|_\mX< M$ for any fixed points of $u=T_\sg(u)$.
\edes

The following lemma establishes Claim 1.

\blemm{Tmaplem} Suppose that there exist $p>s_*/2$ and a constant $M_*$  such that any strong solution $u$ of \mref{gensysfam} satisfies  \beqno{M*defz} \|Du\|_{W^{1,2p}(\Og,\mu)}\le M_*.\eeq Then, there exist  $M,\bg>0$  such that for $\mX=C^{\bg}(\Og)\cap W^{1,2}(\Og)$ the map $T_\sg:\bar{B}_M\to\mX$ is well defined and compact for all $\sg\in[0,1]$. Moreover, $T_\sg$ has no fixed points on $\partial B_M$. \elemm

\bproof 
For some constant $M_0>0$ we consider $u:\Og_{R}\to\RR^m$ satisfying 
\beqno{Xdefstartzz} \|u\|_{C(\Og)}\le M_0,\; \|Du\|_{L^2(\Og)}\le M_0,\eeq
and write the system \mref{Tmapdef} as a linear elliptic system for $w$
\beqno{LSUsys} -\Div(\mathbf{a}(x) Dw) +\mathbf{b}(x)Dw+\mathbf{g}(x)w+Lw =\mathbf{f}(x),\eeq where $ \mathbf{a}(x)=A(x,\sg u)$, $\mathbf{b}(x)=F^{(\sg)}(x,u,Du)$, $\mathbf{g}(x)=f^{(\sg)}(x,u)$, and $\mathbf{f}(x)=\hat{f}(x,0,0)+Lu$.

The matrix $\mathbf{a}(x)$ is then regular elliptic with uniform ellipticity constants by A), AR) because $u$ is bounded. From the theory of {\em linear} elliptic systems it is well known that if the operator $\mathcal{L}(w)=-\Div(\mathbf{a}(x) Dw) +\mathbf{g}(x)w +Lw$ is monotone and there exist positive constants $m$ and $q$ such that  \beqno{LSUcond}\|\mathbf{b}\|_{L^q(\Og)},\; \|\mathbf{g}\|_{L^q(\Og)},\; \|\mathbf{f}\|_{L^q(\Og)}\le m,\; \mbox{$q> n/2$},\eeq then the system \mref{LSUsys}  has a unique weak solution $w$. 

It is easy to find  a matrix $L$ such that $\mathcal{L}(w)$ is monotone. Because the matrix $\mathbf{a}$ is regular elliptic and $\mathbf{g}$ is bounded (see below). We just need to choose a positive definite matrix $L$ satisfying $\myprod{Lw,w}\ge l_0|w|^2$ for some $l_0>0$ and sufficiently large in terms of $M_0$.

Next, we will show that \mref{LSUcond} holds  by F) and \mref{Xdefstartzz}. We consider the two cases  f.1) and  f.2). If f.1) holds then from the definition \mref{Bfdef} there is a constant $C(|u|)$ such that  $$|\mathbf{b}(x)|=|F^{(\sg)}(x,u,\zeta)|\le C(|u|),\; |\mathbf{g}(x)|=|f^{(\sg)}(x,u)|\le C(|u|).$$
 
 From  \mref{Xdefstartzz}, we see that $\|u\|_\infty\le M_0$ and so there is  a constant $m$ depending on $M_0$ such that \mref{LSUcond} holds for any $q$ and $n$.

If f.2) holds then
\beqno{Bfdeff2}|F^{(\sg)}(x,u,\zeta)|\le C(|u|)|\zeta|,\; |f^{(\sg)}(x,u)|\le C(|u|).\eeq Therefore, $\|\mathbf{b}\|_{L^2(\Og)}$ is bounded by $C\|Du\|_{L^2(\Og)}$. Again, if  $n\le3$ then  \mref{Xdefstartzz} implies the condition \mref{LSUcond} for $q=2$.

In both cases,  \mref{LSUsys} (or \mref{Tmapdef}) has a unique weak solution $w$. 
We then define $T_\sg(u)=w$. Moreover, from the regularity theory of linear systems, $w\in C^{\ag_0}(\Og)$ for some $\ag_0>0$ depending on $M_0$.

The bound in the assumption \mref{M*defz}  and \refrem{PSremH} imply that $u$ is H\"older continuous and provide positive constants $\ag, C(M_*)$ such that $\|u\|_{C^{\ag}(\Og)}\le C(M_*)$. Also, the assumption \mref{llgB3} and AR), that $\llg(u),\og$ are bounded from below, yield that $\|Du\|_{L^2(\Og)}\le C(C_0)$. Thus, there is a constant $M_1$, depending on $M_*,C_0$ such that any strong solution $u$ of \mref{gensysfam} satisfies 
\beqno{Xdefstart} \|u\|_{C^{\ag}(\Og)}\le M_1,\; \|Du\|_{L^2(\Og)}\le M_1.\eeq

It is well known that there is a constant $c_0>1$, depending on $\ag$ and the diameter of $\Og$, such that $\|\cdot\|_{C^{\bg}(\Og)}\le c_0\|\cdot\|_{C^{\ag}(\Og)}$ for all $\bg\in(0,\ag)$. We now let $M_0$, the constant in \mref{Xdefstartzz}, be $M=(c_0+1)M_1$ and define $\ag_0$ in the previous argument accordingly. 

Define $\mX=C^{\bg}(\Og)\cap W^{1,2}(\Og)$ for some positive $\bg<\min\{\ag,\ag_0\}$. The space $\mX$ is equipped with the norm $$\|u\|_\mX = \max\{\|u\|_{C^{\bg}(\Og)},\|Du\|_{L^2(\Og)}\}.$$

We now see that $T_\sg$ is well defined and maps the ball $\bar{B}_M$ of $\mX$ into $\mX$. Moreover, from the definition $M=(c_0+1)M_1$, it is clear that $T_\sg$ has no fixed point on the boundary of $B_M$ because such a fixed points $u$ satisfies \mref{Xdefstart} which implies $\|u\|_{\mX}\le c_0M_1<M$. 
 
Finally, we need only show that $T_\sg$ is compact.  If $u$ belongs to a bounded set $K$ of $\bar{B}_M$ then $\|u\|_\mX\le C(K)$ for some constant $C(K)$ and there is a constant $C_1(K)$  such that $\|w\|_{C^{\ag_0}(\Og)}\le C_1(K)$. Thus $T_\sg(K)$ is compact in $C^{\bg}(\Og)$ because $\bg<\ag_0$.  So, we need only show that $T(K)$ is precompact in $W^{1,2}(\Og)$.  We will discuss only the quadratic growth case where \mref{Bfdeff2} holds because the case $\hat{f}$ has linear growth is similar and easier.

First of all, for $u\in K$ we easily see that $\|Dw\|_{L^2(\Og)}$ is uniformly bounded by a constant depending on $K$. The argument is standard by testing the linear system \mref{Tmapdef} by $w$ and using the boundedness of $\|w\|_{L^\infty}$ and $\|u\|_{L^\infty}$, \mref{Bfdeff2}, AR) and Young's inequality.

Let $\{u_n\}$ be a sequence in $K$ and $w_n=T_\sg(u_n)$.  We have, writing $W=w_n-w_m$  $$-\Div(A(x,\sg u_n)DW)=\Div((A(x,\sg u_n)-A(x,\sg u_m)Dw_m)+\Psi_{m,n},$$ where $\Psi_{m,n}$ is defined by $$\barr{lll}\Psi_{m,n}&=&F^{(\sg)}(x,u_n,Du_m)Dw_n-F^{(\sg)}(x,u_m,Du_m)Dw_m+\\&&f^{(\sg)}(x,u_n)u_n-f^{(\sg)}(x,u_m)u_m +L(u_n-u_m+w_m-w_n).\earr$$

Testing the above system with $W$ and using AR), we have
$$\llg_*(K)\mu_*\iidx{\Og}{|DW|^2} \le \iidx{\Og}{[|A(x,\sg u_n)-A(x,\sg u_m)||Dw_m||DW|+|\Psi_{m,n}||W|]}.$$
By Young's inequality, we find a constant $C$ depending on $K$ and $\mu_*$ such that
$$\iidx{\Og}{|DW|^2}\le C\iidx{\Og}{[|A(x,\sg u_n)-A(x,\sg u_m)||Dw_m|^2} +\sup_\Og|W|\|\Psi_{m,n}\|_{L^1(\Og)}.$$
By \mref{Bfdeff2}, it is clear that   $|\Psi_{n,m}|\le C(K)[(|Du_n|+|Du_m|)(|Dw_n|+|Dw_m|)+1]$.
Using the fact that $\|Dw_n\|_{L^2(\Og)}$ and $\|Du_n\|_{L^2(\Og)}$ are uniformly bounded, we see that $\|\Psi_{m,n}\|_{L^1(\Og)}$ is bounded. Hence, $$\iidx{\Og}{|Dw_n-Dw_m|^2}\le C(K)\max\{\sup_\Og|A(x,\sg u_n)-A(x,\sg u_m)|,\sup_\Og|w_n-w_m|\}.$$

Since $u_n,w_n$ are bounded in $C^\bg(\Og)$, passing to subsequences we can assume that $u_n,w_n$ converge in $C^0(\Og)$. Thus, $\|A(x,\sg u_n)-A(x,\sg u_m)\|_\infty,\|w_n-w_m\|_\infty\to0$. We then see from the above estimate that $Dw_n$ converges in $L^{2}(\Og)$. Thus, $T_\sg(K)$ is precompact in $W^{1,2}(\Og)$.

Hence, 
$T_\sg:\mX\to\mX$ is a compact map. The proof is complete. \eproof

We now turn to Claim 2, the hardest part of the proof, and provide a uniform estimate for the fixed points of $T_\sg$ and justify the key assumption \mref{M*defz} of \reflemm{Tmaplem}. To this end, we first have the following lemma. 

\blemm{Dulocbound} A fixed point of $T_\sg$ is also a strong solution of \mref{gensysfam}. 
 \elemm
 
The  lemma can easily follow from the results in \cite{Gius}. However, there are some remarks need be made here as we are also considering quadratic growth $\hat{f}$ in this paper and our system is not necessarily variational.

\brem{Dulocboundrem} If $u$ is a fixed point of $T_\sg$ in $\mX$ then it solves \mref{gensysfam} weakly and is continuous. Thus, $u$ is bounded and belongs to $VMO(\Og)$. By AR), the system \mref{gensysfam} is regular elliptic. We can adapt the proof in \cite{GiaS}, which deals with parabolic systems, to our elliptic case.
If $\hat{f}$ satisfies a quadratic growth in $Du$ then,  because $u$ is bounded, the condition \cite[(0.4)]{GiaS} that $|\hat{f}|\le a|Du|^2+b$ is satisfied here. The proof of \cite[Theorems 2.1 and 3.2]{GiaS} assumed the 'smallness condition' (see \cite[(0.6)]{GiaS}) $2aM<\llg_0$, where $M=\sup |u|$. This 'smallness condition' was needed because only weak bounded solutions, which are not necessarily continuous, were considered in \cite{GiaS}. In our case,  $u$ is continuous so that we do not require this 'smallness condition'.  Indeed, a careful checking of the arguments of the proof in \cite[Lemma 2.1 and page 445]{GiaS} shows that if $R$ is small and one knows that the solution $u$ is continuous then these argument still hold as long we can absorb the integrals involving $|Du|^2, |Dw|^2$ (see the estimate after \cite[(3.7)]{GiaS}) on the right hand sides to the left right hand sides of the estimates. Thus, \cite[Theorems 2.1 and 3.2]{GiaS} apply to our case and yield that $u\in C^{a}(\Og)$ for all $a\in(0,1)$ and that, since $A(x,u)$ is differentiable,  $Du$ is locally H\"older continuous in $\Og$. Therefore, $u$ is also a strong solution (see \cite[Chapter 10]{Gius}).\erem

Thanks to \reflemm{Dulocbound}, we need only consider a strong solution $u$ of \mref{gensysfam} and establish \mref{M*defz} for some $p>s_*/2$. Because the data of \mref{Tmapdef} satisfy the structural conditions A), F) with the same set of constants and the assumptions of the theorem are assumed to be uniform for all $\sg\in[0,1]$, we will only present the proof for the case $\sg=1$ in the sequel.

Let $u$ be a strong solution  of \mref{e1} on $\Og$.  We begin with an energy estimate for $Du$. For $p\ge1$ and any ball $B_s$ with center $x_0\in\bar{\Og}$ we denote $\Og_s=B_s\cap\Og$ and 
\beqno{Hdef}\ccH_{p}(s):=
\iidmu{\Og_s}{\llg(u)|Du|^{2p-2}|D^2u|^2},\eeq
\beqno{Bdef}\ccB_{p}(s):= \iidmu{\Og_s}{\frac{|\llg_u(u)|^2}{\llg(u)}|Du|^{2p+2}},\eeq \beqno{Cdef}\ccC_{p}(s):=\iidmu{\Og_s}{(|f_u(u)|+\llg(u))|Du|^{2p}},\eeq and \beqno{cIdef}\mccIee_{\og,p}(s):=\iidx{B_s}{(\llg(u)|Du|^{2p}|D\og_0|^2+|f(u)||Du|^{2p-1}|D\og_0|\og_0)}.\eeq

\blemm{dleANSenergy}  Assume A), F). 
Let $u$ be any strong solution  of \mref{gensys} on $\Og$ and  $p$ be any number in $[1,\max\{1,s_*/2\}]$.

There is a constant $C$, which depends only on the parameters in A) and F), such that for any two concentric balls $B_s,B_t$ with center $x_0\in\bar{\Og}$ and $s<t$ \beqno{keydupANSenergy}\ccH_{p}(s)\le C\ccB_{p}(t)+C(1+(t-s)^{-2})[\ccC_{p}(t)+\mccIee_{\og,p}(t)].\eeq

\elemm

\bproof The proof is similar to the energy estimate of $Du$ for the parabolic case in \cite[Lemma 3.2]{dleANS}. As we consider $A, \hat{f}$ depending on $x$, there are some extra terms in the estimate which need some extra attention.
We follow the proof of  \cite[Lemma 3.2]{dleANS} with $W=U=u$ and $\bg(u)=1$. Roughly speaking, we differentiated the system in $x$ to obtain
\beqno{ga2zzz} -\Div(A(x,u)D^2u+A_u(x,u)DuDu+A_x(x,u)Du)=D\hat{f}(x,u,Du).\eeq

For any two concentric balls $B_s,B_t$, with $s<t$, let $\psi$ be a cutoff function for $B_s,B_{t}$. That is, $\psi$ is a $C^1$ function satisfying $\psi\equiv1$ in $B_s$ and $\psi\equiv0$ outside $B_t$ and $|D\psi|\le 1/(t-s)$. We then
test \mref{ga2zzz} with $|Du|^{2p-2}Du\psi^2(x)$ and obtain, using integration by parts and Young's inequality
\beqno{keyenergystart}\barr{ll}\lefteqn{
\iidx{\Og_t}{\llg(u)|Du|^{2p-2}|D^2u|^2\psi^2\og}\le  C\iidx{\Og_t}{[\frac{|\llg_u(u)|^2}{\llg(u)}|Du|^{2p+2}+|D\psi|^{2}\llg(u)|Du|^{2p}]\og}}\hspace{1cm}&\\& +\iidx{\Og_t}{[|A_x(x,u)||Du|^{2p-1}|D^2u|\psi^2+C(p)|D\hat{f}(x,u,Du)||Du|^{2p-1}\psi^2]}.\earr\eeq

Here, integrals in the first line of \mref{keyenergystart} result from the same argument in the proof of \cite[Lemma 3.2]{dleANS} using the spectral gap condition SG) we are assuming here (see also \reflemm{SGrem} in the Appendix).  Meanwhile, the integral in the last line comes from the presence of $x$ in $A,\hat{f}$ and we will discuss it below. 

We consider the integrand $|A_x(x,u)||Du|^{2p-1}|D^2u|\psi^2$ on the right hand side of \mref{keyenergystart}. As $|A_x(x,u)|\le C_*\llg(u)|D\og|$ and $\og=\og_0^2$, we have to deal with the integral $$C_*\iidx{\Og}{\llg(u)|Du|^{2p-1}|D^2u|\psi^2|D\og_0|\og_0}.$$
By Young's inequality this integral can be estimated by   \beqno{Axest2} \frac12\iidx{\Og}{ \llg(u)|Du|^{2p-2}|D^2u|^2\psi^2\og_0^2}+\frac12C_*^2 \iidx{\Og}{\llg(u)|Du|^{2p}\psi^2|D\og_0|^2}.\eeq

Next, we consider the integral of $|D\hat{f}(x,u,Du)||Du|^{2p-1}\psi^2$. First, if $\hat{f}$ has a linear growth in $Du$ then by f.1) in F) with $p=Du$, $$|D\hat{f}(x,u,Du)| \lesssim  (\llg(u)|D^2u|+ |\llg_u(u)||Du|^2)\og+\llg(u)|Du||D\og|+ |f_u(u)||Du|\og +|f(u)||D\og|.$$ Therefore, using   Young's inequality,  we get
\beqno{f1duest}\barr{ll}\lefteqn{|D\hat{f}(x,u,Du)||Du|^{2p-1}\lesssim\llg(u)|Du|^{2p}|D\og_0|^2+|f(u)||Du|^{2p-1}|D\og_0|\og_0+}\hspace{.2cm}&\\& \left\{\eg\llg(u)|Du|^{2p-2}|D^2u|^2+ C(\eg)\llg(u)|Du|^{2p}+\frac{|\llg_u(u)|^2}{\llg(u)}|Du|^{2p+2}+|f_u(u)||Du|^{2p}\right\}\og_0^2.\earr\eeq

Similarly, if $\hat{f}$ has a quadratic growth in $Du$ then by f.2) in F) $$\barr{lll}|D\hat{f}(x,u,Du)| &\lesssim&  (|\llg_u(u)||Du||D^2u|+ |\llg_{uu}(u)||Du|^3)\og+ |\llg_u||Du|^2|D\og|+\\&&
|f_u(u)||Du|\og +|f(u)||D\og|.\earr$$

We then have to deal with three extra terms which can be handled by  Young's inequality and the assumption \mref{llgquadcond} that $|\llg_{uu}(u)|\llg(u)\lesssim |\llg_u(u)|^2$ as follows.
$$|\llg_u(u)||Du||D^2u||Du|^{2p-1}\le \eg\llg(u)|Du|^{2p-2}|D^2u|^2 + C(\eg)\frac{|\llg_u(u)|^2}{\llg(u)}|Du|^{2p+2},$$
$$|\llg_u||Du|^2|D\og||Du|^{2p-1}\lesssim \frac{|\llg_u(u)|^2}{\llg(u)}|Du|^{2p+2}\og+\llg(u)|Du|^{2p}|D\og|^2\og^{-1},$$
$$|\llg_{uu}(u)||Du|^3|Du|^{2p-1}\lesssim \frac{|\llg_u(u)|^2}{\llg(u)}|Du|^{2p+2}.$$ We then get the same terms as in \mref{f1duest} for the linear growth case. 

Finally, we use the definitions \mref{Hdef}-\mref{cIdef} to see that the integral in the last line of \mref{keyenergystart} can be estimated by
\beqno{energyDu}(\eg+\frac12)\iidmu{\Og_t}{\llg(u)|Du|^{2p-2}|D^2u|^2\psi^2}+C(\eg)[\ccB_p(t)+(1+(t-s)^{-2})(\ccC_p(t)+\mccIee_{\og,p}(t))].\eeq
We then choose $\eg$  sufficiently small  so that the first integral can be absorbed to the left hand side of \mref{keyenergystart}. We then obtain \mref{keydupANSenergy} and complete the proof. \eproof

Next, if the condition \mref{Keymu0} of M.1) holds then we combine the energy estimate and the local Gagliardo-Nirenberg inequality \mref{GNlocog11mcoro} to have the following stronger estimate.

\blemm{dleANSprop}  In addition to the assumptions of \reflemm{dleANSenergy}, we suppose that M.1) holds for some $p$. That is, for any given  
$\mu_0>0$ there exist a constant $C_0$ and a positive $R_{\mu_0}$ sufficiently small in terms of the constants in A) and F) such that
\beqno{Keymu01} \sup_{x_0\in\bar{\Og}}[\myPi_p^{\ag}]_{\bg+1,\Og_R(x_0)}\le C_0,\;\|K(u)\|_{BMO(\Og_{R}(x_0),\mu)}^2 \le \mu_0.\eeq

Then for sufficiently small $\mu_0$ there is a constant $C$ depending only on the parameters of A) and F) such that for $2R<R_{\mu_0}$ we have \beqno{keydupANSppb}\ccB_{p}(R)+\ccH_{p}(R)\le C(1+R^{-2})[\ccC_{p}(2R)+\mccIee_{\og,p}(2R)].\eeq

 \elemm

\bproof  Recall the energy estimate \mref{keydupANSenergy} in \reflemm{dleANSenergy} 
\beqno{keyiteration1} \ccH_{p}(s)\le C\ccB_{p}(t)+C(1+(t-s)^{-2})[\ccC_{p}(t)+\mccIee_{\og,p}(t)], \; 0<s<t.\eeq 

We apply \refcoro{GNlocalog1mcoro} to estimate $\ccB_{p}(t)$, the integral  on the right hand side of \mref{keyiteration1}. We let $\LLg(u)=\llg^\frac{1}{2}(u)$ in \refcoro{GNlocalog1mcoro} and note that $\myPi_p$ defined there is now comparable to the $\myPi_p=\llg^{p+\frac12}(u)|\llg_u(u)|^{-p}$ in M.1). We compare the definitions \mref{Idefmt} and \mref{I0*mt} with those in \mref{Hdef}-\mref{Cdef} to see that for $U(x)=u(x)$ $$\ccB_{p}(t)=I_1(t,x_0),\;\ccC_{p}(t)=I_0(t,x_0),\;\ccH_{p}(t)=I_2(t,x_0).$$
Hence, for any $\eg>0$ we can use \mref{GNlocog11mcoro} obtain a constant $C$ such that (using the bound $[\myPi_p^{\ag}]_{\bg+1,B_{R_{\mu_0}}(x_0)\cap\Og}\le C_0$ and the definitions of $\mu_0$ in \mref{Keymu01}  and $C(\eg,U,\myPi)$  in \refcoro{GNlocalog1mcoro})
$$ \ccB_{p}(s) \le \eg\ccB_{p}(t)+C\eg^{-1}\mu_0\ccH_{p}(t)+C\eg^{-1}\mu_0(t-s)^{-2}\ccC_{p}(t)\quad 0<s<t\le R_{\mu_0}.$$

Define $F(t):=\ccB_{p}(t)$, $G(t):=\ccH_{p}(t)$, $g(t):=\ccC_{p}(t)$ and $\eg_0=\eg+C\eg^{-1}\mu_0$. The above yields \beqno{keyiteration2}F(s)\le \eg_0[F(t)+G(t)] +C(t-s)^{-2}g(t).\eeq 

Now, for $h(t):=\mccIee_{\og,p}(t)$ the energy estimate \mref{keyiteration1}  implies
 \beqno{Giter}G(s)\le C[F(t)+(1+(t-s)^{-2})(g(t)+h(t))].\eeq

As $\eg_0=\eg+C\eg^{-1}\mu_0$, it is clear that we can choose and fix some $\eg$ sufficiently small and then $\mu_0$ small in terms of $C,\eg$ to have $2C\eg_0<1$. Thus, if $\mu_0$ is sufficiently small in terms of the  constants in A),F), then we can apply a simple iteration argument \cite[Lemma 3.11]{dleANS} to the two inequalities \mref{keyiteration2} and \mref{Giter} and obtain for $0<s<t\le R_{\mu_0}$ $$ F(s)+G(s) \le C(1+(t-s)^{-2})[g(t)+h(t)].$$ 

For any $R<R_{{\mu_0}}/2$ we take $t=2R$ and $s=\frac32R$ in the above to obtain $$\ccB_{p}(\frac32R)+\ccH_{p}(\frac32R)\le C(1+R^{-2})[\ccC_{p}(2R)+\mccIee_{\og,p}(2R)].$$ Combining this and \mref{keyiteration1} with $s=R$ and $t=\frac32R$, we see that
$$\ccB_{p}(R)+\ccH_{p}(R)\le C(1+R^{-2})[\ccC_{p}(2R)+\mccIee_{\og,p}(2R)].$$
This is \mref{keydupANSppb} and the proof is complete. \eproof

Finally, we have the following lemma giving a uniform bound for strong solutions.

\blemm{dleANSpropt0}  Assume as in \reflemm{dleANSprop} and AR). 
We assume further that there exists a constant $C_0$ such that for $r_*={\pi_*}/({\pi_*}-2)$ and some $r_0>r_*$ \beqno{llgmainhyp} \|\llg^{-1}(u)\|_{L^{r_*}(\Og,\mu)},\;\||f_u(u)|\llg^{-1}(u)\|_{L^{r_0}(\Og,\mu)}\le C_0 ,\eeq
\beqno{hypoiterpis1}\iidmu{\Og}{(|f_u(u)|+\llg(u))|Du|^{2}}\le C_0,\eeq and
\beqno{hypoiterpis1z}\iidmu{\Og}{|f_u(u)|(1+f(u)^{s_*}|f_u(u)|^{-s_*})}\le C_0.\eeq

Then there exist $p>s_*/2$ and a constant $M_*$ depending only on the parameters of A) and F), $\mu_0$, $R_{\mu_0}$, $C_0$ and the geometry of $\Og$  such that 
\beqno{keydupANSt0}\iidmu{\Og}{|Du|^{2p}}\le   M_*.\eeq
 \elemm

\bproof First of all, by the condition AR), there is a constant $C_\og$ such that $|D\og_0|\le C_\og\og_0$ and therefore we have from the the definition \mref{cIdef} that \beqno{cIdef1}\mccIee_{\og,p}(s)\le C_\og\iidx{B_s}{(\llg(u)|Du|^{2p}+f(u)|Du|^{2p-1})\og_0^2}.\eeq By Young's inequality, $f(u)|Du|^{2p-1}\lesssim |f_u(u)||Du|^{2p}+(f(u)|f_u(u)|^{-1})^{2p}|f_u(u)|$. It follows easily from the assumption  \mref{hypoiterpis1z} that the integral of $(f(u)|f_u(u)|^{-1})^{2p}|f_u(u)|$ is bounded by $C_0$ for any $p\in[1,s_*/2]$. We then have from \mref{keydupANSppb} and \mref{cIdef1} that
\beqno{keydupANSppbbb}\ccB_{p}(R)+\ccH_{p}(R)\le C(1+R^{-2})[\ccC_{p}(2R)+C_0].\eeq

The main idea of the proof is to show that the above estimate is self-improving in the sense that if it is true for some exponent $p\ge1$ then it is also true for $\cg_*p$ with some fixed $\cg_*>1$ and $R$ being replaced by $R/2$. To this end, assume that for some $p\ge1$ we can find a constant $C(C_0,R,p)$ such that \beqno{piter1} \ccC_{p}(2R)\le C(C_0,R,p).\eeq 

It then clearly follows from \mref{keydupANSppbbb}, and the definition of $\ccB_{p}(R),\ccH_{p}(R)$ that \beqno{keydupANSpp1} \iidmu{\Og_R}{[V^2+
|DV|^2]}\le C(C_0,R,p), \quad \mbox{where $V=\llg^\frac12(u)|Du|^{p}$}.\eeq

Let $\pi_*$ be the exponent in the Poincar\'e-Sobolev inequality \mref{PSineq}. We have 
\beqno{genPSpi}\iidmu{\Og_R}{|V|^{\pi_*}}\lesssim  \iidmu{\Og_R}{|DV|^{2}}+\mypar{\iidmu{\Og_R}{|V|^2}}^{{\pi_*}/2}.\eeq We see that \mref{keydupANSpp1} and the above inequality imply
\beqno{piter1az} \iidmu{\Og_R}{\llg^{{\pi_*}/2}(u)|Du|^{p {\pi_*}}} \le C(C_0,R,p).\eeq

Let $\cg_*\in(1,{\pi_*}/2)$.  We  denote $\bg_*=\pi_*/(\pi_*-2\cg_*)$ and $g(u)=\max\{|f_u(u)|,\llg(u)\}$. We write $g(u)|Du|^{2\cg_* p}= \llg(u)^{\cg_*}|Du|^{2\cg_* p}g(u)\llg(u)^{-\cg_*}$ and  and use H\"older's inequality and \mref{piter1az} to get
\beqno{piter1z1}\iidmu{\Og_R}{g(u)|Du|^{2\cg_* p}}\le C(C_0,R,p)^\frac{2\cg_*}{\pi_*}\mypar{\iidmu{\Og_R}{(g(u)\llg(u)^{-\cg_*})^{\bg_*}}}^{1-\frac{2\cg_*}{{\pi_*}}}.\eeq Again, as $(g(u)\llg(u)^{-\cg_*})^{\bg_*}=(g(u)\llg(u)^{-1})^{\bg_*}\llg(u)^{-(\cg_*-1)\bg_*}$, the last integral can be bounded via H\"older's inequality by $$\mypar{\iidmu{\Og_R}{(g(u)\llg(u)^{-1})^{\bg_*\ag}}}^\frac1\ag\mypar{\iidmu{\Og_R}{(\llg(u)^{-(\cg_*-1)\bg_*\ag'}}}^\frac{1}{\ag'}$$
By the assumption \mref{llgmainhyp},  $|f_u(u)|\llg^{-1}(u)$ is in $L^r_0(\Og,\mu)$ for some $r_0>r_*=\pi_*/(\pi_*-2)$ and $\llg^{-1}(u)\in L^{r_*}(\Og,\mu)$.  We can find $\ag,\cg_*>1$ such that $\ag\bg_*<r_0$ and $(\cg_*-1)\bg_*\ag'<r_*$ and see that, from the definition $g(u)=\max\{|f_u(u)|,\llg(u)\}$, the above integrals are bounded by some constant $C(C_0)$. Hence, this fact and \mref{piter1z1} imply\beqno{piter1z}\ccC_{\cg_*p}(R)=\iidmu{\Og_R}{(|f_u(u)|+\llg(u))|Du|^{2\cg_* p}}\le C(C_0,R,p).\eeq

We just show that if \mref{piter1} holds true for some $p\ge 1$ and $R>0$ then \mref{piter1z}  provides some fixed  $\cg_{*}>1$ such that \mref{piter1} remains true for the new exponent $\cg_{*} p$  and $R/2$. By the assumption \mref{hypoiterpis1}, \mref{piter1} holds for $p=1$. It is now clear that, as long as the energy estimate \mref{keydupANSenergy} is valid by \reflemm{dleANSenergy}), we can repeat the argument $k_0$ times to find  a number $p>s_*/2$ such that \mref{piter1} holds. It follows  that 
there is a constant $C$ depending only on the parameters of A) and F), $\mu_0$, $R_{\mu_0}$ and $k_0$ 
such that for some $p>s_*/2$
\beqno{keydupANSt0z}\iidmu{\Og_{R_0}}{\llg^{{\pi_*}/2}(u)|Du|^{{\pi_*}p}}\le   C\mbox{ for  $R_0=2^{-k_0}R_{\mu_0}$}.\eeq

We now write $|Du|^{2p}=\llg(u)|Du|^{2p}\llg^{-1}(u)$ and have
$$\iidmu{\Og_{R_0}}{|Du|^{2p}}\le \mypar{\iidmu{\Og_{R_0}}{\llg^{{\pi_*}/2}(u)|Du|^{{\pi_*}p}}}^\frac{2}{{\pi_*}}\mypar{\iidmu{\Og_{R_0}}{\llg(u)^{-{{\pi_*}/({\pi_*}-2)}}}}^{1-\frac{2}{{\pi_*}}}.$$ By \mref{llgmainhyp},  the last integral in the above estimate is bounded.
Using \mref{keydupANSt0z} and summing the above inequalities over a finite covering of balls $B_{R_0}$ for $\Og$, we find a constant $C$, depending also on the geometry of $\Og$, and obtain the desired estimate \mref{keydupANSt0}. The lemma is proved. \eproof

\brem{hypoiterpis1zrem} The condition \mref{hypoiterpis1z} is void if $\hat{f}$ is independent of $x$. Indeed, it was used only to estimate $\mccIee_{\og,p}$, which results from \mref{Giter} in the proof of \reflemm{dleANSprop}, and obtain \mref{keydupANSppbbb}. If $\hat{f}$ is independent of $x$ then $\mccIee_{\og,p}=0$ from the proof of the energy estimate in \reflemm{dleANSenergy} (see \mref{f1duest}) so that \mref{hypoiterpis1z} is not needed. \erem

\brem{murem1} It is also important to note that the estimate of \reflemm{dleANSpropt0}, based on those in  \reflemm{dleANSenergy}, \reflemm{dleANSprop}, is  {\em independent} of lower/upper bounds of the function $\llg_*$ in AR) but the integrals in M.0). The assumption AR) was used only in \reflemm{Tmaplem} to define the map $T_\sg$ and \reflemm{Dulocbound} to show that fixed points of $T_\sg$ are strong solutions. \erem

We are ready to provide the proof of the main theorem of this section. 

{\bf Proof of \reftheo{gentheo1}:} 
It is now clear that the assumptions M.0) and M.1) of our theorem allow us to
apply \reflemm{dleANSpropt0} and obtain a priori uniform bound  for any continuous strong solution $u$ of \mref{gensysfam}. The uniform estimate \mref{keydupANSt0} shows that the assumption \mref{M*defz} of \reflemm{Tmaplem} holds true so that the map $T_\sg$ is well defined and compact on a ball $\bar{B}_M$ of $\mX$ for some $M$ depending on the bound $M_*$ provided by \reflemm{dleANSpropt0}. Combining with \reflemm{Dulocbound}, the fixed points of $T_\sg$ are strong solutions of the system \mref{gensysfam} so that $T_\sg$ does not have a fixed point on the boundary of $\bar{B}_M$. Thus, by the Leray-Schauder fixed point theorem, $T_1$ has a fixed point in $B_M$ which is a strong solution to \mref{gensys}. The proof is complete.
\eproof

{\bf Proof of \refcoro{gentheo1coro}:} We just need to show that \reflemm{dleANSpropt0} remains true with the condition \mref{hypoiterpis1z}, which is \mref{hypoiterpis1z00}, being replaced by \mref{hypoiterpis1z11} and \mref{hypoiterpis1z11z}. We revisit the proof of \reflemm{dleANSpropt0}. First of all,  we use Young's inequality in the estimate \mref{cIdef1} for $\mccIee_{\og,p}(s)$ to get $f(u)|Du|^{2p-1}\lesssim |f_u(u)||Du|^{2p}+(f(u)|f_u(u)|^{-1})^{2p}|f_u(u)|$. Using \mref{specfcond} and \mref{hypoiterpis1z11},  \mref{keydupANSppbbb} now yields $$\ccB_{p}(R)+\ccH_{p}(R)\le C(1+R^{-2})[\ccC_{p}(2R)+C_0+I_p(2R)], \quad I_p(s)=\iidmu{\Og_s}{|u|^{2p}|f_u(u)|}.$$ 

As in \mref{piter1}, using the same idea, we will first assume for some $p\ge1$ that \beqno{piter111} \ccC_{p}(2R)+I_p(2R)\le C(C_0,R,p)\eeq and show that this assumption is self-improving, i.e., if it holds for some $p\ge1$ then it remains true for $\cg_*p$ for some fixed $\cg_*>1$ and $2R$ being replaced by $R$. We need only consider $I_p(R)$ and assume first that \beqno{improvep1} \iidmu{\Og_R}{(|f_u(u)|+\llg(u))|u|^{2}}\le C(C_0,R).\eeq

Denote $g(u)=|f_u(u)|$ and $\cg_*=2-\frac{1}{r_0}-\frac{2}{\pi_*}$. Since $r_0>\frac{\pi_*}{\pi_*-2}$ it is clear that $\cg_*>1$. We then define $r_1=\frac{\pi_*}{(\cg_*-1)(\pi_*-2)}$, $r_2=\frac{\pi_*}{2\cg_*}$ and note that
$$\frac{1}{r_0}+\frac{1}{r_1}+\frac{1}{r_2}=\frac{1}{r_0}+\frac{(\cg_*-1)(\pi_*-2)}{\pi_*}+\frac{2\cg_*}{\pi_*}=\frac{1}{r_0}+\cg_*-1+\frac{2}{\pi_*}=1.$$ Therefore, writing
$g(u)|u|^{2p\cg_*}=g(u)\llg^{-1}(u)\llg^{1-\cg_*}(u)\llg^{\cg_*}|u|^{2p\cg_*}$, we can use  H\"older's inequality to see that $I_{p\cg_*}(R)$ can be bounded by \beqno{fullgup}\mypar{\iidmu{\Og_R}{(g(u)\llg(u)^{-1})^{r_0}}}^\frac{1}{r_0}\mypar{\iidmu{\Og_R}{\llg(u)^{-\frac{\pi_*}{\pi_*-2}}}}^\frac{1}{r_1}\mypar{\iidmu{\Og_R}{\llg(u)^{\frac{\pi_*}{2}}|u|^{p\pi_*}}}^\frac{1}{r_2}.
\eeq
Thanks to \mref{llgmainhyp}, the first two integrals are bounded by $C(C_0)$. We estimate the third integral. Let $U=\llg(u)^{\frac{\pi_*}{4}}|u|^\frac{p\pi_*}{2}$. Because $|\llg_u(u)||u|\lesssim\llg(u)$, $|DU|^2\lesssim \llg(u)^{\frac{\pi_*}{2}}|u|^{p\pi_*-1}|Du|$. By Poincar\'e and Young's inequalities, we have for any $\bg\in(0,1]$
$$\iidmu{\Og_R}{U^2}\le R^2\iidmu{\Og_R}{\llg(u)^{\frac{\pi_*}{2}}(\eg|u|^{p\pi_*}+C(\eg)|Du|^{p\pi_*})}+R^{n-n/\bg}\mypar{\iidmu{\Og_R}{U^{2\bg}}}^\frac{1}{\bg}.$$

Choosing $\eg$ small and $\bg=2/\pi_*$, we see that the right hand side is bounded by $$ CR^2\iidmu{\Og_R}{\llg(u)^{\frac{\pi_*}{2}}|Du|^{p\pi_*}}+L_p(R)^\frac{\pi_*}{2},\; \mbox{where } L_p(s)=\iidmu{\Og_s}{\llg(u)|u|^{2p}}.$$

Putting these estimates together and using \mref{piter1az}, which holds because of \mref{piter111}, we see that the third integral in \mref{fullgup} is bounded by a constant $C(C_0,R,p)$ if $L_p(R)\le C(C_0,R,p)$. We then have $I_{\cg_* p}(R)\le C(C_0,R,p)$.

On the other hand, we can show that the estimate $L_p(R)\le C(C_0,R,p)$ is also self improving. We repeat the argument in \mref{fullgup} with $g(u)=\llg(u)$ and, of course, $r_0=\infty$ to see that if $L_p(R)\le C(C_0,R,p)$ then $L_{\cg_* p}(R)\le C(C_0,R,p)$ for $\cg_*=2-\frac{2}{\pi_*}>1$ because $\pi_*>2$.

Hence, the estimate \mref{piter111} remains true for $p,2R$ being replaced by $\cg_*p,R$ respectively for some fixed $\cg_*>1$. The proof of \reflemm{dleANSpropt0} can continue with \mref{improvep1} replacing \mref{hypoiterpis1z}.

Finally, we show that \mref{improvep1} comes from the assumptions \mref{hypoiterpis1} and \mref{hypoiterpis1z11}. From \mref{fullgup} with $p=1$, we see that $I_1(R)$ can be estimated in terms of the integral of $\llg(u)^{\frac{\pi_*}{2}}|u|^{\pi_*}$ over $\Og_R$. The latter can be estimated by $L_1(R)$ via the Poincar\'e-Sobolev inequality \mref{genPSpi}, using $V=\llg(u)^{\frac{1}{2}}|u|$ and the assumption on the integrability of $|DV|^2\sim\llg(u)|Du|^2$ in \mref{hypoiterpis1}. Hence, we need only consider $L_1(R)$. By the interpolation inequality, we have $$\iidmu{\Og_R}{\llg(u)|u|^{2}}\le R^2\iidmu{\Og_R}{\llg(u)|Du|^2}+R^{n-n/\bg}\mypar{\iidmu{\Og_R}{(\llg(u)|u|^{2})^\bg}}^\frac{1}{\bg}.$$ We then use \mref{hypoiterpis1} and \mref{hypoiterpis1z11} to see that $L_1(R)\le C(C_0,R)$ and conclude the proof. \eproof

\section{Proof of the main results} \label{proofmainsec} \eqnoset

We now present the proof of the results in \refsec{res} which are in fact just applications of the main technical theorem with different choices of the map $K$ verifying the conditions H) and M.1). Again, since the systems in \mref{e1famzzz} satisfy the same set of conditions uniformly with respect to $\sg\in[0,1]$ we will only present the proof for $\sg=1$ in the sequel.

{\bf Proof of \reftheo{dleNL-mainthm}:} This theorem is a consequence of \refcoro{gentheo1coro} whose integrability conditions in M.0) are already assumed in \mref{llgmainhyp009} and \mref{llgmainhyp0099}. We need only verify the condition M.1) of \refcoro{gentheo1coro} (or \reftheo{gentheo1}) with the map $K$ being defined by
\beqno{Kmapdef} K_{\eg_0}(u)=(|\log(|U|)|+\eg_0)|U|^{-1}U, \quad U=[\llg_0+|u_i|]_{i=1}^m.\eeq This map satisfies for any $\eg_0>0$ the condition H) of \reftheo{gentheo1} (see \reflemm{H2lemm} and \refrem{H2lemmrem} in the Appendix). We need only check the condition M.1). Because $$\|K_{\eg_0}(u)\|_{BMO(B_R,\mu)}\le \|K_0(u)\|_{BMO(B_R,\mu)}+\eg_0\||U|^{-1}U\|_{BMO(B_R,\mu)},$$ and $\||U|^{-1}U\|_{BMO(B_R,\mu)}=1$ and so $\|K_{\eg_0}(u)\|_{BMO(B_R,\mu)}<\mu_0$ for any given $\mu_0>0$ if $\eg_0<\mu_0/2$ and $R$ is small, thanks to the assumption \mref{mainloghyp}. Therefore, the smallness condition \mref{Keymu0} of M.1) holds.

Next, we consider $\myPi_p=\llg^{p+\frac12}(u)|\llg_u(u)|^{-p}\sim (\llg_0+|u|)^{k_p}$ for some number $k_p$ because $\llg(u), |\llg_u(u)|$ have polynomial growths.  As we assume in \mref{mainloghyp} that $\log(\llg_0+|u|)$ has small BMO norm in small balls,  $[\myPi_p^\ag]_{\bg+1,B_R}\sim [(\llg_0+|u|)^{k_p\ag}]_{\bg+1,B_R}$ is bounded (see \reflemm{Westlemm} following this proof) for any given $\ag,\bg>0$ if $R$ is sufficiently small.  Thus, M.1) is verified and \refcoro{gentheo1coro} applies here to complete the proof. \eproof

\newc{\mmc}{\mathbf{c}}

We have the following lemma which was used in the above proof to establish that $[\myPi_p^\ag]_{\bg+1,B_R}$ is bounded. This lemma will be frequently referred to in the rest of this section. 

\blemm{Westlemm} Let $\mu$ be a doubling measure and $U$ be a nonnegative function on a ball $B\subset \Og$. There is a constant $\mmc_2$ depending only on the doubling constant of $\mu$ such that for any given $l$ and $\bg>0$ if $[\log U]_{*,\mu}$ is sufficiently small then $[U^l]_{\bg+1}\le \mmc_2^{1+\bg}$.
\elemm

\bproof We first recall the John-Nirenberg inequality (\cite[Chapter 9]{Graf}): For any BMO($\mu$) function $v$ 
there are constants $\mmc_1,\mmc_2$, which depend only on the doubling constant of $\mu$, such that  
\beqno{JNineq}\mitmu{B}{e^{\frac{\mmc_1}{[v]_{*,\mu}}|v-v_B|}} \le \mmc_2.\eeq

For any $\bg>0$ we know that $e^v$ is an $A_{\bg+1}$ weight with $[e^v]_{\bg+1}\le \mmc_2^{1+\bg}$ (e.g. see \cite[Chapter 9]{Graf}) if  \beqno{Acgcond}\sup_B \mitmu{B}{e^{(v-v_B)}} \le \mmc_2,\; \sup_B \mitmu{B}{e^{-\frac{1}{\bg}(v-v_B)}} \le \mmc_2.\eeq

It is clear that \mref{Acgcond} follows from \mref{JNineq} if $\mmc_1[v]_{*,\mu}^{-1}\ge \max\{1,\bg^{-1}\}$. Therefore, for $v= l\log U$ we see that if $[\log U]_{*,\mu}\le \mmc_1\bg|l|^{-1}$ then $[U^l]_{\bg+1}\le \mmc_2^{1+\bg}$. Hence, for any given $l$ and $\bg>0$ if $[\log U]_{*,\mu}$ is sufficiently small then $[U^l]_{\bg+1}\le \mmc_2^{1+\bg}$.
\eproof

{\bf Proof of \refcoro{dlec1coro}:}  We just need to show that the assumption \mref{pis20az5} implies \mref{mainloghyp} of \reftheo{dleNL-mainthm}. For $U=[\llg_0+|u_i|]_{i=1}^m$ we have from the calculation in \refrem{H2lemmrem} that
$$(K_0)_u(u) = \frac{|\log(|U|)|}{|U|}\left[I + \left(\frac{\mbox{sign}(\log(|U|))}{|\log(|U|)|}-1\right)\zeta\zeta^T\right]\mbox{diag}[\mbox{sign} (u_i)],$$ where $\zeta=|U|^{-1}U$. It is clear that $|(K_0)_u(u)|\le C(1+|\log(|U|))|U|^{-1}$. Since $|U|$ is bounded from below by $\llg_0$, for any $\ag\in(0,1)$ there is a constant $C(\ag)$ such that $|(K_0)_u(u)|\le C(\ag)|U|^{-\ag}$. Therefore, $|D(K_0(u))|$ and $|D(\log(\llg_0+|u|))|$ can be bounded by $C(\llg_0+|u|)^{-\ag}|Du|$. 
It follows from the assumption \mref{pis20az5} that \beqno{pis20az5zz} \iidmu{\Og}{|DK_0(u)|^{n}},\; \iidmu{\Og}{|D(\log(\llg_0+|u|))|^{n}}\le C(C_0).\eeq

From the Poincar\'e-Sobolev inequality, using the assumption that $\og$ is bounded from below and above by AR) 
$$\left(\mitmu{B_r}{|K_0(u)-K_0(u)_{B_r}|^2}\right)^\frac12 \le C(n)
\left(\iidmu{B_r}{|D(K_0(u))|^{n}}\right)^\frac{1}{n}.$$ The  continuity of the integral of $|D(K_0(u))|^{n}$ and the uniform bound \mref{pis20az5zz} show that the last integral is small if $r$ is. The same argument applies to the function $\log(\llg_0+|u|)$. We then see that the BMO norms of $K_0(u)$ and  $\log(\llg_0+|u|)$ are small in small balls, and so \mref{mainloghyp} of \reftheo{dleNL-mainthm} holds. The proof is complete.
\eproof

For the proof of \refcoro{2dthm1} we first need the following lemma.

\blemm{2dlemma} Suppose A), F) and \mref{FUDU112z} if $\hat{f}$ has a quadratic growth in $Du$ with $\eg_0$ being sufficiently small. For any $s$ satisfying \beqno{sconds}s>-1 \mbox{ and } C_*^{-1}>s/(s+2)\eeq and any ${\ag_0}\in(0,1)$ we have for $U:=\llg_0+|u|$ that \beqno{nequal2est}\iidmu{\Og}{U^{k+s}|Du|^2} \le C\mypar{\iidmu{\Og}{U^{{\ag_0}(k+s+2)}}}^\frac1{\ag_0}+C\iidmu{\Og}{U^sf(u)|u|}.\eeq \elemm

\bproof  
Let $X=[\llg_0+|u_i|]_{i=1}^m$ and test the system with $|X|^su$ to get \beqno{llgDu2}\iidmu{\Og}{\myprod{A(u)Du,D(|X|^su)}} \le \iidmu{\Og}{\myprod{\hat{f}(u,Du),|X|^su}}.\eeq  We note that $\myprod{A(u)Du,D(|X|^su)}=\myprod{A(u) DX,D(|X|^sX)}$ so that, by  the assumption \mref{sconds} on $s$, there is $c_0>0$ such that $\myprod{A(u)Du,D(|X|^su)}\ge c_0\llg(u)|X|^s|DX|^2$ (see \mref{SGrem0} in the Appendix). Because $|DX|=|Du|$ and $|X|\sim U$, the above yields \beqno{fuduu}\iidmu{\Og}{\llg(u)U^s|Du|^2}\le C\iidmu{\Og}{U^s\myprod{\hat{f}(u,Du),u}}.\eeq

If $\hat{f}$ satifies f.1) then a simple use of Young's inequality gives $$|\myprod{\hat{f}(u,Du),u}| \le \eg\llg(u)|Du|^2 + C(\eg)\llg(u)|u|^2+f(u)|u|.$$ 

If f.2) holds with \mref{FUDU112z} then $|\myprod{\hat{f}(u,Du),u}| \le C|\llg_u(u)||u||Du|^2 + f(u)|u|$. Because $|\llg_u(u)||u|\lesssim \llg(u)$, we obtain the above inequality again with $\eg=C\eg_0$. Therefore, if $\eg_0$ is sufficiently small then \mref{fuduu} and the fact that $\llg(u)\sim U^k$  imply $$\iidmu{\Og}{U^{k+s}|Du|^2} \le C\iidmu{\Og}{U^{k+s+2}}+C\iidmu{\Og}{U^sf(u)|u|}.$$
We apply the interpolation inequality $\|w\|_{L^2(\Og,\mu)}^2 \le \eg\|Dw\|_{L^2(\Og,\mu)}^2+C(\eg)\|w^{\ag_0}\|_{L^2(\Og,\mu)}^{1/{\ag_0}}$ with ${\ag_0}\in(0,1)$ and $w=U^{(k+s+2)/2}$ to the first integral on the right hand side, noting also that $|Dw|^2\sim U^{k+s}|Du|^2$. For $\eg$ sufficiently small, we derive  \mref{nequal2est} from the above estimate  and complete the proof. \eproof

{\bf Proof of \refcoro{2dthm1}:} We apply \refcoro{dlec1coro} here. We will verify first the condition \mref{pis20az5} and then the integrability assumptions \mref{llgmainhyp009} and \mref{llgmainhyp0099}.

In the sequel we denote $U=\llg_0+|u|$. From \mref{nequal2est} of \reflemm{2dlemma} we see that  if  there exist positive numbers $\ag_0\in(0,1)$, $C_0>0$ and $s$ satisfies \mref{sconds} such that if
\beqno{nequal2conda}\|U^{p}\|_{L^1(\Og)}\le C_0,\;p=\ag_0(s+k+2)\mbox{ and }p=s+l+1,\eeq then \mref{nequal2est}, with the assumption that $f(u)\lesssim U^l$, implies  \beqno{nequal2condb}\iidx{\Og}{U^{k+s}|Du|^2}\le C(C_0).\eeq

We will first show that  the assumption \mref{keyn2lnorm0}, that $\|u\|_{L^{l_0}(\Og,\mu)}\le C_0$ for some $l_0>\max\{l,l-k-1\}$, implies \mref{nequal2condb} for some $s=s_0>\max\{-1,-k-2\}$. Indeed, for any such $l_0$ we have $s_0+l+1\le l_0$ if $s_0$ is close to $\max\{-1,-k-2\}$. Moreover, $s_0$ satisfies \mref{sconds}. This clearly holds if $s_0\le0$, i.e. $k\ge-2$. Otherwise, the assumption $k>-2C_*/(C_*-1)$ yields that $C_*^{-1}>s_0/(s_0+2)$ if $s_0$ is close to $-k-2>0$.  Thus, \mref{nequal2est} holds for such $s_0$. We also choose $\ag_0\in(0,1)$ sufficiently small such that $\ag_0(s_0+k+2)\le l_0$. With these choices of $\ag_0, s_0$ and the assumption \mref{keyn2lnorm0}, we see that \mref{nequal2conda} and then \mref{nequal2condb} hold for $s=s_0$.

As $k+s_0>-2$, we can find $\ag\in(0,1)$ such that $-2\ag\le k+s_0$. Therefore, \mref{nequal2condb} with $s=s_0$ yields \mref{pis20az5} of \refcoro{dlec1coro} for $n=2$ because $$\iidx{\Og}{(\llg_0+|u|)^{-2\ag}|Du|^2} \le \iidx{\Og}{(\llg_0+|u|)^{k+s_0}|Du|^2}\le C(C_0).$$

We now check the integrability conditions \mref{llgmainhyp009} and \mref{llgmainhyp0099} of \reftheo{dleNL-mainthm} which read
\beqno{llgmainhyp00} \|\llg^{-1}(u)\|_{L^{\frac n2}(\Og,\mu)},\;\||f_u(u)|\llg^{-1}(u)\|_{L^{r_0}(\Og,\mu)}\le C_0,\eeq
\beqno{llgB30}\iidmu{\Og}{(|f_u(u)|+\llg(u))|Du|^{2}}\le C_0,\eeq 
\beqno{llgB40}\iidmu{\Og}{(\llg(u)|u|^2)^{\bg_0}}\le C_0.\eeq

Because $n=2$, we have the inequality $\|w\|_{L^q(\Og)} \le C\|Dw\|_{L^2(\Og)}+C\|w^\bg\|_{L^1(\Og)}^\frac{1}{\bg}$ which holds for all $q\ge1$ and $\bg\in(0,1)$. Applying this to $w=|U|^{(k+s_0)/2+1}$ and using \mref{nequal2condb} and the assumption \mref{keyn2lnorm0}, we see that $\|U^q\|_{L^1(\Og)}\le C(C_0)$ for all $q\ge1$. By H\"older's inequality this is also true for $q\ge0$. It is also true for $q<0$ because $U\ge\llg_0>0$. We then have\beqno{n2allq} \|U^q\|_{L^1(\Og)}\le C(C_0,\llg_0,q) \quad \mbox{for all $q$}.\eeq

The above then immediately implies the integrability conditions \mref{llgmainhyp00} and \mref{llgB40} because $\llg(u)$ and $|f_u(u)|$ are powers of $U$.

Similarly, \mref{n2allq} implies that \mref{nequal2conda} holds for any $p$ so that \mref{nequal2condb} is valid if $s\ge0$ and $C_*^{-1}>s/(s+2)$.  To verify  \mref{llgB30} we need to find a constant $C(C_0)$ such that \beqno{llgB30zz}\iidx{\Og}{(\llg_0+|u|)^{l-1}|Du|^2}+\iidx{\Og}{(\llg_0+|u|)^k|Du|^{2}}\le C(C_0).\eeq Letting $s=0$ in \mref{nequal2condb}, we see that the second integral on the left hand side is bounded by a constant $C(C_0)$. If $l\le 1+k$ then the first integral in \mref{llgB30zz} is bounded by the second one and we obtain the desired bound. If $l>k+1$ we let $s=l-k-1$ in \mref{nequal2condb}. The condition on $s$ in \mref{sconds} holds because $$\frac{s}{s+2}<C_*^{-1} \Leftrightarrow \frac{l-k-1}{l-k+1}<C_*^{-1} \Leftrightarrow l-k<\frac{C_*+1}{C_*-1},$$ which is assumed in \mref{klcond}. Hence, \mref{llgB30zz} holds. We have verified all assumptions of \reftheo{dleNL-mainthm} and the  proof is complete. \eproof

We now give the proof of \reftheo{n3SKT}. The case $n=2$ is similar and easier so that we will consider only $n=3,4$. The proof is again based on \reftheo{gentheo1} using the new map  
\beqno{Kn34def}K(u)=[K_i(u)]_{i=1}^m, \mbox{ where } K_i(u)=\log(\llg_0^{k+1}+|P_i(u)|).\eeq

We first check the conditions H) of \reftheo{gentheo1} in the following lemma.
\blemm{newKlem} For any $k\ge-1$ there exists a constant $C(k)$ such that for $\mathbb{K}(u)=K_u^{-1}(u)^T$ $$|\mathbb{K}(u)|\le C(k) \llg(u)|\llg_u(u)|^{-1},\; |\mathbb{K}_u(u)|\le C(k).$$
\elemm
\bproof As $A(u)=P_u(u)$, we have from the definition \mref{Kn34def} that  $$K_u(u)=[\frac{\partial K_i(u)}{\partial u_j}(K_i(u))^{-1}]=\mbox{diag}[\mbox{sign}P_i(u)(\llg_0^{k+1}+|P_i(u)|)^{-1}]A(u),$$ and so 
$K_u^{-1}(u)=A(u)^{-1}\mbox{diag}[\llg_0^{k+1}\mbox{sign}P_i(u)+P_i(u)]$. We  show that there is a cosntant $C(k)$ such that $|K_u^{-1}(u)|$ is bounded by $C(k)\llg(u)|\llg_u(u)|^{-1}$. Since  $\llg(u), C_*\llg(u)$ are the smallest and largest eigenvalues of $A(u)$, $|A^{-1}(u)|\sim \llg^{-1}(u)$. Using the facts that $\llg(u)\sim (\llg_0+|u|)^k$, $|\llg_u(u)|\sim (\llg_0+|u|)^{k-1}$ and $(\llg_0+|u|)^{k+1}\ge \llg_0^{k+1}$ (because $k\ge-1$),  it follows easily that $\llg_0^{k+1}|\llg_u(u)|\lesssim \llg(u)^{2}$. We collect these facts in the following. \beqno{Allgk1}|P_i(u)|\lesssim |A(u)||u|,\;|A^{-1}(u)|\sim \llg^{-1}(u),\; \llg_0^{k+1}|\llg_u(u)|\lesssim \llg(u)^{2}.\eeq
We then have  $|A(u)^{-1}\mbox{diag}[P_i(u)]|\lesssim |u|\lesssim \llg(u)|\llg_u(u)|^{-1}$ and $\llg_0^{k+1}|A(u)^{-1}|\sim \llg_0^{k+1}\llg(u)^{-1}\lesssim \llg(u)|\llg_u(u)|^{-1}$. Therefore  $|\mathbb{K}(u)|=|K_u^{-1}(u)^T|\lesssim \llg(u)|\llg_u(u)|^{-1}$. 

In addition, as $|\mathbb{K}_u(u)|= |(K_u^{-1}(u))_u|$ and $A(u)=\mbox{diag}[(P_i(u))_u]$ $$|\mathbb{K}_u(u)|\lesssim |A^{-1}(u)|^2|A_u(u)||\mbox{diag}[\llg_0^{k+1}+|P_i(u)|]|+1.$$
Using \mref{Allgk1} and the fact that $|A_u(u)|\lesssim |\llg_u(u)|$,  we easily see that $|\mathbb{K}_u(u)|$ is bounded by some constant $C(k)$. The lemma is proved. \eproof

Next, we need to show that $K(u)$ has small BMO norm in small balls by establishing that $DK(u)\in L^n(\Og)$ and using the Poincar\'e-Sobolev inequality as in the proof of \refcoro{dlec1coro}.  To this end, we need
the following lemma.
\blemm{DU6lem} Assume that there exist $r_0>n/2$ and $C_0$ such that
\beqno{llgn3}  \|f_u(u)\llg^{-1}(u)\|_{L^{r_0}(\Og)}\le C_0.\eeq
Then for any $\bg_0\in(0,1]$  there exists a constant $C(C_0,\bg_0)$,  such that \beqno{DUUest}\|DP(u)\|_{L^\frac{2n}{n-2}(\Og)}\le C(C_0,\bg_0)\||P(u)|^{\bg_0}\|_{L^1(\Og)}.\eeq \elemm
\bproof In the sequel, we write $U=[U_i]_{i=1}^m$, $U_i=P_i(u)$. Multiplying the $i$-th equation in \mref{genSKT} by $-\Delta U_i$, integrating over $\Og$ and summing the results, we get
$$\iidx{\Og}{|\Delta U|^2}=-\sum_i\iidx{\Og}{\myprod{B_i(u,Du),\Delta U_i}}-\sum_i\iidx{\Og}{\myprod{f_i(u),\Delta U_i}}.$$ Applying integration by parts to the last integral, we have $$\iidx{\Og}{|\Delta U|^2}=-\sum_i\iidx{\Og}{\myprod{B_i(u,Du),\Delta U_i}}+\sum_{i,j}\iidx{\Og}{\myprod{(f_i)_{u_j}(u)Du_j,D U_i}}.$$

The condition A) implies $\llg(u)|Du|^2\le \myprod{A(u)Du,Du}=\myprod{DU,Du}$ and so Young's inequality yields $\llg(u)|Du|^2\le \frac12\llg^{-1}(u)|DU|^2+\frac12\llg(u)|Du|^2$. We then have $\llg(u)|Du|\le |DU|$.
Using this fact, the assumption that $|B_i(u,Du)|\le C\llg(u)|Du|$ and applying Young's inequality to the first integral on the right hand side of the above, we get $$\|\Delta U\|_{L^2(\Og)}^2\lesssim\iidx{\Og}{|DU|^2}+\iidx{\Og}{|f_u(u)\llg^{-1}(u)||DU|^2}.$$
By H\"older's inequality and \mref{llgn3}, the last integral is estimated by
$$ \mypar{\iidx{\Og}{|f_u(u)\llg^{-1}(u)|^{r_0}}}^\frac{1}{r_0}\|DU\|_{L^{2r'_0}(\Og)}^2 \le C(C_0)\|DU\|_{L^{2r'_0}(\Og)}^2.$$

Using Schauder's inequality $\|D^2 U\|_{L^2(\Og)}\le C\|\Delta U\|_{L^2(\Og)}$, we obtain from the above two inequalities that
\beqno{D2UPP}\|D^2U\|_{L^2(\Og)}^2\lesssim \|DU\|_{L^2(\Og)}^2+C(C_0)\|DU\|_{L^{2r'_0}(\Og)}^2.\eeq

We recall the following interpolating Sobolev inequality: for any  $\eg>0$  \beqno{intineq}\|W\|_{L^p(\Og)}\le \eg\|DW\|_{L^2(\Og)} + C(\eg)\|W^\bg\|_{L^1(\Og)}^\frac{1}{\bg} \mbox{ for any $p\in[1,\frac{2n}{n-2})$ and $\bg\in(0,1]$}.\eeq 

Because $r_0>n/2$, $2r'_0<2n/(n-2)$ so that we can apply \mref{intineq} to $W=DU$ with $p=2$, $p=2r'_0$, $\bg=1$ and $\eg$ is sufficiently small in \mref{D2UPP} to see that  $\|D^2U\|_{L^2(\Og)}^2\le C(C_0)\|DU\|_{L^1(\Og)}^2$. As $\|D U\|_{L^1(\Og)}\le \eg\|D^2U\|_{L^2(\Og)} + C(\eg)\|U\|_{L^1(\Og)}$,
we obtain for small $\eg$ that $\|D^2U\|_{L^2(\Og)}\le C(C_0)\|U\|_{L^1(\Og)}$.
Sobolev's embedding theorem then yields $$\|D U\|_{L^\frac{2n}{n-2}(\Og)}\le C(C_0)\|U\|_{L^1(\Og)}.$$
Applying \mref{intineq} again, with $W=|U|$, $p=1$ and $\bg=\bg_0$, to estimate the norm $\|U\|_{L^1(\Og)}$ and noting that $\|DU\|_{L^2(\Og)}\lesssim\|DU\|_{L^\frac{2n}{n-2}(\Og)}$, we obtain \mref{DUUest}. 
\eproof

{\bf Proof of \reftheo{n3SKT}:} The proof is again based on \reftheo{gentheo1}. The assumptions i) and ii) state \beqno{uintcondk}\left\{\barr{ll}\|u\|_{L^1(\Og)}\le C_0 & \mbox{if $k\ge0$,}\\\|u\|_{L^{-kn/2}(\Og)}\le C_0 &\mbox{if $k\in[-1,0)$,}\earr \right.\eeq and clearly imply
 \beqno{fuhyp0z}\|\llg^{-1}(u)\|_{L^{\frac{n}{2}}(\Og)}\le C(C_0,\llg_0).\eeq This and the assumption \mref{fuhyp0} provide the integrability condition \mref{llgmainhyp0} of M.0). Concerning the integrability condition \mref{llgB3} in M.0), we use H\"older's inequality, writing $\llg(u)|Du|^2= \llg^{-1}(u)\llg^2(u)|Du|^2$, and \mref{fuhyp0z} and \mref{DUUest} to have
 $$\iidx{\Og}{\llg(u)|Du|^2}\le\|\llg^{-1}(u)\|_{L^\frac{n}{2}(\Og)}\|\llg(u)Du\|_{L^\frac{2n}{n-2}(\Og)}^2\le C_0\|DP(u)\|_{L^\frac{2n}{n-2}(\Og)}^2\le C(C_0).$$
 
 Similarly, the integral of $|f_u(u)||Du|^2$ can be estimated by
 $$\|\llg^{-2}(u)f_u(u)\|_{L^\frac{n}{2}(\Og)}\|\llg(u)Du\|_{L^\frac{2n}{n-2}(\Og)}^2\le  C(C_0)\|\llg^{-2}(u)f_u(u)\|_{L^\frac{n}{2}(\Og)}.$$
 
 If $k\ge0$ then $\llg^{-2}(u)|f_u(u)|\le C(\llg_0)\llg^{-1}(u)|f_u(u)|$ so that the last norm in the above is bounded, thanks to \mref{fuhyp0}. If $k<0$ then this norm is bounded by the assumption \mref{fuhyp0k}. We conclude that the condition \mref{llgB3} in M.0) holds.

We discuss the condition M.1). First of all, 
because $|P(u)|\le \llg(u)|u|\le (\llg_0+|u|)^{k+1}$, \mref{uintcondk} also shows that for any positive and sufficiently small $\bg_0$ \beqno{Puhypz} \|(\llg_0^{k+1}+|P(u)|)^{\bg_0}\|_{L^1(\Og)}\le C(\llg_0,C_0).\eeq

Next, using the Poincar\'e-Sobolev inequality as in the proof of \refcoro{dlec1coro}, we show that $K(u)$ has small BMO norm in small balls by estimating
\beqno{3dlog}\|DK_i(u)\|_{L^n(\Og)} \le \|(\llg_0^{k+1}+|P_i(u)|)^{-1}DP_i(u)\|_{L^n(\Og)}\le C(\llg_0)\|DP_i(u)\|_{L^n(\Og)}.\eeq
Because $n\le4$, $n\le 2n/(n-2)$. By the assumption \mref{fuhyp0}, $\|f_u(u)\llg^{-1}(u)\|_{L^{r_0}(\Og)}\le C_0$, \reflemm{DU6lem} shows that $\|DP(u)\|_{L^n(\Og)}\le C(C_0,\bg_0)\||P(u)|^{\bg_0}\|_{L^1(\Og)}$. This and \mref{Puhypz} and \mref{3dlog} provide a constant $C(\llg_0,C_0)$ such that  $\|DK(u)\|_{L^n(\Og)}\le C(\llg_0,C_0)$. We then see that $K(u)$ has small BMO norm in small balls.

Concerning the weight $\myPi_p$, because $|Du|\lesssim\llg^{-1}(u)|DP(u)|$, $k+1\ge0$, we have
$$|D(\log(\llg_0+|u|))|=\frac{|Du|}{\llg_0+|u|}\lesssim \frac{|DP(u)|}{(\llg_0+|u|)\llg(u)}\sim \frac{|DP(u)|}{(\llg_0+|u|)^{k+1}}\le \llg_0^{-k-1}|DP(u)|.$$ \reflemm{DU6lem} then shows that $D\log(\llg_0+|u|)\in L^n(\Og)$ so that $\log(\llg_0+|u|)$ is has small BMO norm in small balls. \reflemm{Westlemm} applies to yield that $[\myPi_p^\ag]_{\bg+1,B_R}\sim [(\llg_0+|u|)^{k_p\ag}]_{\cg,B_R}$ is bounded for any given $\ag,\bg>0$ if $R$ is sufficiently small. The assumption M.1) of \reftheo{gentheo1} is verified.

We thus establish the conditions of \reftheo{gentheo1} and complete the proof. \eproof

\section{Appendix}\label{appsec}\eqnoset

Let $m,l$ be any integers. For $X=[X_i]_{i=1}^m$, $X_i\in \RR^l$ and for any $C^1$ function $k:\RR^+\to\RR$ we consider the maps
\beqno{kdef} K(X)=k(|X|)|X|^{-1}X,\;\zeta=|X|^{-1}X=[\zeta_i]_{i=1}^m.\eeq

We see that $D_X(|X|)=\zeta$ and $D_X\zeta=|X|^{-1}(I-\zeta\zeta^T)$, where $\zeta\zeta^T=[\myprod{\zeta_i,\zeta_j}]$. Hence, 
$$D_XK(X)= k(|X|)D_X\zeta+D_Xk(|X|)\zeta^T=k(|X|)|X|^{-1}(I-\zeta\zeta^T)+k'(|X|)\zeta\zeta^T.$$ 
We then introduce the notations \beqno{kudef}b=k'(|X|)|X|/k(|X|),\; \ccK(\zeta)=I+(b-1)\zeta\zeta^T.\eeq Therefore, the calculation for $D_XK(X)$ yields \beqno{kxdef} D_XK(X)=k(|X|)|X|^{-1}(I+(b-1)\zeta\zeta^T)=k(|X|)|X|^{-1}\ccK(\zeta).\eeq

If $k(|X|)\ne 0$ and $k'(|X|)\ne 0$ then $\ccK(\zeta)$ is invertible. 
We can use the Serman-Morrison formula $(I+wv^T)^{-1}=I-(1+v^Tw)^{-1}wv^T$, setting $w=(b-1)\zeta$ and $v=\zeta$, to see that  \beqno{DXKinv} (D_XK(X))^{-1}=\frac{|X|}{k(|X|)}(I+(b^{-1}-1)\zeta\zeta^T).\eeq 

Otherwise, if $k(|X|)=0$ (resp. $k'(|X|)=0$) then $D_XK(X)= k'(|X|)\zeta\zeta^T$ (resp. $D_XK(X)= k(|X|)|X|^{-1}(I-\zeta\zeta^T)$) and $D_XK(X)$ is not invertible.

The following lemma was used in the checking of the condition H) for the map $K_{\eg_0}(u)$ in the proof of \reftheo{dleNL-mainthm}.
\blemm{H2lemm} For any $\eg_0,\llg_0>0$ let $k(t)=|\log(t)|+\eg_0$ and $X(u)=[\llg_0+|u_i|]_{i=1}^m$ in \mref{kdef}. There exists a constant $C(\eg_0)$ such that the map $\mathbb{K}(u)=(K_u(X(u))^{-1})^T$ satisfies \beqno{H2check}|\mathbb{K}(u)|\le C(\eg_0)|X|,\; \|\mathbb{K}_u(u)\|_{L^\infty(\RR^m)}\le C(\eg_0).\eeq

\elemm

\bproof As $k'(t)=\mbox{sign}(\log t) t^{-1}$, we have $$b^{-1}=k(|X|)(k'(|X|)|X|)^{-1}=\mbox{sign}(\log(|X|)) (|\log(|X|)|+\eg_0).$$ 

Define $X_u=\mbox{diag}[\mbox{sign} u_i]$.  We have from \mref{DXKinv} and the definition $X=[\llg_0+|u_i|]_{i=1}^m$ that
$$\mathbb{K}(u)=(X_u)^{-1}\frac{|X|}{|\log(|X|)|+\eg_0}(I+(\mbox{sign}(\log(|X|)) (|\log(|X|)|+\eg_0)-1)\zeta\zeta^T).$$

As $\eg_0>0$, we easily see that $|\mathbb{K}(u)|\le C(\eg_0)|X|$ for some constant $C(\eg_0)$. A straightforward calculation also shows that $\|\mathbb{K}_u(u)\|_{L^\infty(\RR^m)}\le C(\eg_0)$.  \eproof

\brem{H2lemmrem}
If $\llg(u)\sim (\llg_0+|u|)^k\sim |X(u)|^k$, with $k\ne0$, and $\LLg(u)= \llg^\frac12(u)$ and $\Fg(u)=|\LLg_u(u)|$. We then have $\LLg(u)\Fg^{-1}(u)\sim |X(u)|$ and $\Fg(u)|\Fg_u(u)|^{-1}\sim |X(u)|$. We obtain from \mref{H2check} that $|\mathbb{K}(u)|\lesssim \LLg(u)\Fg^{-1}(u)$ and  $|\mathbb{K}(u)||\Fg_u(u)|\lesssim\Fg(u)$.
Thus, the assumptions on the map $K$ for the local Gagliardo-Nirenberg inequality are verified here. \erem

We then need that $K(U)$ is BMO and $\myPi^\ag$ is a weight.  By \mref{kudef} $$K_U(U) = \frac{|\log(|X|)|+\eg_0}{|X|}\left(I + (\frac{X_s}{|\log(|X|)|+\eg_0}-1)\zeta\zeta^T\right)X_U, \quad X=[\llg_0+|U_i|]_{i=1}^m.$$

Recall that $$\myPi=\LLg^{p+1}(U)\Fg^{-p}(U)\sim |U|^p\llg^\frac12(U),\; \ag>2/(p+2),\; \bg<p/(p+2).$$

We note some properties of $K$.

\blemm{Klemm} Let $X_i:\Og\to \RR^l$ be $C^1$ maps on a domain $\Og\subset \RR^n$ for $i=1,\ldots,m$. If $b=k'(|X|)|X|/k(|X|)>0$ then \beqno{DXDK}\myprod{DX,D(K(X))}\ge0.\eeq Moreover, for $\ag=1-(\frac{b-1}{b+1})^2$
\beqno{DXDK1}\myprod{DX,D(K(X))}\ge \ag^\frac12 |DX||D(K(X))|. \eeq

\elemm

\bproof By \mref{kxdef} $D(K(X))=\mmk(|X|)\ccK(\zeta)DX$, where $\mmk(t):=k(t)t^{-1}$, and so
$$\myprod{DX,D(K(X))}=\mmk(|X|)(|DX|^2+(b-1)\myprod{DX,\zeta\zeta^T DX}).$$ Note that $|\zeta\zeta^T|\le1$ so that $\myprod{DX,D(K(X))}\ge0$ if $s=b-1>-1$. This gives \mref{DXDK}.

Since $\zeta\zeta^T$ is a projection, i.e. $(\zeta\zeta^T)^2=\zeta\zeta^T$, we have, setting $J=\myprod{\zeta\zeta^T DX,DX}$
$$\barr{lll}|D(K(X))|^2&=&|K_X(X)DX|^2=\mmk^2(|X|)\myprod{\ccK(\zeta)DX,\ccK(\zeta)DX}\\&=&
\mmk^2(|X|)(|DX|^2+(2s+s^2)J).\earr$$
Hence, we can write $\myprod{DX,D(K(X))}^2-\ag |DX|^2|D(K(X))|^2$ as
$$\mmk^2(|X|)[(1-\ag)|DX|^4+(2s-\ag(2s+s^2))|DX|^2 J+s^2J^2].$$
If we choose $\ag=1-(\frac{s}{s+2})^2$ then the above is $\mmk^2(|X|)\left(\frac{s}{s+2}|DX|^2-sJ\right)^2\ge0$.
Therefore $\myprod{DX,D(K(X))}^2\ge \ag |DX|^2|D(K(X))|^2$. This and \mref{DXDK} yield \mref{DXDK1}. \eproof

We now consider a matrix $A$ satisfying for some positive $\llg,\LLg$ and any vector $\chi$\beqno{Aellcond} \myprod{A\chi,\chi}\ge \llg|\chi|^2,\quad |A\chi|\le \LLg|\chi|.\eeq 

\blemm{AuK} Assume \mref{Aellcond} and that $b>0$ (see \mref{kudef}). For $\kappa=\llg/\LLg^2$ and $\nu=\llg/\LLg$ we have $$\myprod{\kappa A DX,D(K(X))}\ge (\ag^\frac12-(1-\nu^2)^\frac12) |DX||D(K(X))|.$$
\elemm

\bproof From \mref{Aellcond} with $\chi=DX$, we note that $$\barr{lll}|\kappa ADX-DX|^2 &=& \kappa^2|ADX|^2-2\kappa\myprod{ADX,DX}+|DX|^2\\&\le& (\kappa^2\LLg^2-2\kappa\llg+1)|DX|^2=(1-\nu^2)|DX|^2.\earr$$

Therefore, using \mref{DXDK1}$$\barr{lll}\myprod{\kappa ADX,D(K(X))}&=&\myprod{\kappa ADX-DX,D(K(X))}+\myprod{DX,D(K(X))}\\&\ge& -|\kappa ADX-DX||DK(X)|+\ag^\frac12 |DX||D(K(X))|.\earr$$ As $|\kappa ADX-DX|^2\le (1-\nu^2)|DX|^2$, we obtain the lemma. \eproof

Let $k(t)=|t|^{s+1}$ then $K(X)=|X|^sX$ and $b=s+1$. The above lemma then gives the following result which was used in the energy estimate.
\blemm{SGrem} Assume \mref{Aellcond}. If $s>-1$ and $\nu=\frac{\llg}{\LLg}>\frac{s}{s+2}$, then  \beqno{SGrem0}\myprod{A DX,D(|X|^sX)}\ge  c_0\frac{\LLg^2}{\llg}|X|^s|DX|^{2},\eeq where $c_0=(1-(\frac{s}{s+2})^2)^\frac12-(1-\nu^2)^\frac12>0$.
\elemm

\bibliographystyle{plain}

\end{document}